\documentclass[oneside]{amsart}

\usepackage{geometry}
\usepackage{appendix}
\usepackage{amsmath}
\usepackage{graphicx}
\usepackage{lineno}
\usepackage{array}
\usepackage{longtable}
\usepackage{natbib}
\usepackage{blindtext}
\usepackage{float}
\usepackage{amssymb}
\usepackage{makecell}
\usepackage{amssymb}
\usepackage{bm}
\usepackage{url}
\usepackage{xcolor}
\usepackage[ruled, lined, linesnumbered, commentsnumbered, longend]{algorithm2e}
\def\BibTeX{{\rm B\kern-.05em{\sc i\kern-.025em b}\kern-.08em
    T\kern-.1667em\lower.7ex\hbox{E}\kern-.125emX}}

\theoremstyle{definition}

\theoremstyle{remark}

\numberwithin{equation}{section}

\SetCommentSty{mycommfont} 
\title[Non-parametric clustering of SDE] {Statistical learning of nonlinear stochastic differential equations from 
non-stationary time series using variational clustering}

\author{Vyacheslav Boyko$^1$}
\address[Vyacheslav Boyko]{$^1$Department of Mathematics and Computer Sciences\\
        Freie Universitaet Berlin\\
        Germany, Arnimallee 6-14195 Berlin}
\email[V. Boyko]{vyacheslav.boyko@fu-berlin.de}

\author{Sebastian Krumscheid$^2$}
\address[Sebastian Krumscheid]{$^2$Faculty of Mathematics, Computer Science and Natural Sciences\\
	RWTH Aachen University\\
	Germany, Templergraben 55, 52062 Aachen}
\email[S. Krumscheid]{krumscheid@uq.rwth-aachen.de}

\author{Nikki Vercauteren$^{3}$}
\address[Nikki Vercauteren]{$^3$Department of Geosciences\\
	 University of Oslo\\
	Norway, P.O. Box 1072 Blindern, 0316 Oslo}
\email[N. Vercauteren]{nikki.vercauteren@geo.uio.no}
\date{\today}

\begin{document}
\maketitle
\begin{abstract}
Parameter estimation for non-stationary stochastic differential equations (SDE)
with an arbitrary nonlinear drift, and nonlinear diffusion is accomplished in
combination with a non-parametric clustering methodology.  Such a model-based
clustering approach includes a quadratic programming (QP) problem with equality
and inequality constraints. We couple the QP problem to a closed-form likelihood
function approach based on suitable Hermite expansion to approximate the
parameter values of the SDE model. The classification problem provides a smooth
indicator function, which enables us to recover the underlying temporal
parameter modulation of the one-dimensional SDE. The numerical examples show
that the clustering approach recovers a hidden functional relationship between
the SDE model parameters and an additional auxiliary process. The study builds
upon this functional relationship to develop closed-form, non-stationary,
data-driven stochastic models for multiscale dynamical systems in real-world
applications.

\end{abstract} \section{Introduction}\label{intro}
Clustering approaches are reliable tools for the analysis of complex time series
with the eventual goal to extract patterns and thereby gain an understanding of
the processes in natural and physical science \citep{zhao_data_2005,
	ayenew_hierarchical_2009, maione_research_2019}. If one supplements the
classification method with a model structure, the clustering becomes model-based
\citep{liao_clustering_2005}. This feature enhances the modeling potential and
provides a way for out-of-sample predictions. The mathematical model takes
different structures depending on the assumptions. One valuable requirement is a
formulation that allows an analytical study after the model has been identified.
Therefore, the estimation of a continuous-type model is preferable because the
intermediate discretization step entails further post-processing effort. An
additional advantage of the continuous model lies in further developing and
embedding the continuous-type estimates into a multiscale system. This is an
arduous specification since the clustering method should be general enough to
handle a wide range of dynamics, including nonstationarities and nonlinearity.

In a series of works \cite{horenko_finite_2010, metzner_analysis_2012,
	pospisil_scalable_2018} introduced and developed an efficient non-parametric
model-based clustering framework, which proved to be a successful analysis tool
in atmospheric sciences \citep{horenko_identification_2010, okane_changes_2013,
	vercauteren_clustering_2015, franzke_systematic_2015, risbey_metastability_2015,
	okane_dynamics_2016, vercauteren_statistical_2019, boyko_multiscale_2020}
molecular dynamics \citep{gerber_inference_2014} and computational finance
\citep{putzig_optimal_2010}. The paradigm of the approach is based on two
assumptions. The non-stationary time series is assumed to be represented by a
non-autonomous model with time-dependent parameters. The timescales of the
parameters are assumed to be much longer than the fluctuation timescale of the
time series. In analogy, if we consider the oscillations of the model parameters
and the time series in the frequency domain, then we would assume a significant
scale separation.

Essentially, the clustering method solves two inverse problems: one to determine
the parameters of a mathematical model and one to classify the
cluster-respective parameter sets. The number of parameter sets, i.e., clusters,
is a hyperparameter, which resolves the parameter variability.  The
classification is interpreted in terms of the location in time, simultaneously
providing the periods and their optimal, locally-stationary parameter sets. The
two inverse problems are coupled through the so-called fitness function, which
allows combining both problems in a single minimization procedure. Such a
variational framework is based on a regularized non-convex clustering algorithm,
which allows efficient non-parametric modeling of non-stationary time series.

The main difference among existing variational frameworks similar to the one
introduced by \cite{horenko_finite_2010}  lies in the assumed underlying
mathematical model. One approach is based on the vector autoregressive factor
model with exogenous variables \citep{horenko_identification_2010}, where the
discrete-time linear model defines the model structure. Conversely, the approach
of \citep{metzner_analysis_2012} uses a Markov-regression model, while a Markov
chain Monte Carlo method is used by \cite{de_wiljes_adaptive_2013} . The
variational framework has also been combined with the probability
distribution-based modeling by utilizing the generalized extreme value
distribution \citep{kaiser_inference_2014}. A more general distribution-based
clustering is constructed based on the maximum entropy principle. This framework
provides a memoryless, multidimensional, and non-stationary modeling
\citep{horenko_computationally_2020}. The framework is computationally scalable
and enables clustering of both temporal and spatial parameter modulations
\citep{kaiser_data-based_2015}. Such efficiency is achieved by the suitable
utilization of an appropriate finite-dimensional function space (e.g., finite
element spaces) to approximate the solution of the classification problem.

Similar to the approaches mentioned above, we adopt the variational framework to
combine it with the model structure defined by stochastic differential equations
(SDE) in a time-continuous formulation. The novelty of our approach is that we
have an analytical, closed-form approximation of the likelihood function that we
use as the model-observation misfit (or fitness) function; this is in contrast
to the Euclidean norm used by \citep{horenko_identification_2010}. The maximum
likelihood estimator (MLE) is then the result of minimizing the misfit (or
maximizing the fitness). In this way, we avoid an explicit discretization step
of the mathematical SDE model. Additionally, we incorporate the functional
modeling of the noise term in the classification algorithm, improving the
accuracy of the classification-modeling task. In particular, it allows to model
nonlinear noise processes and not just additive ones. The MLE is constructed by
using the transition density of a general, one-dimensional nonlinear SDE. The
closed-form expression of the transition density is approximated via the Hermite
polynomials following \cite{ait-sahalia_maximum_2002}. For example, this
approach showed to be successful in modeling complex movement patterns of marine
predators and climate transitions \citep{krumscheid_data-driven_2015}.  Here we
do not deal with the problem of determining an appropriated functional form of
the underlying SDE, instead of with the question: how can one arrive from a
clustering approach to a self-contained predictive model?

The goal is to present a data-driven stochastic modeling strategy for
parameterizing unresolved degrees of freedom in a complex system. We consider a
multiscale system composed of a model representing the macroscale, while a
low-dimensional model represents the fast degrees of freedom. Theoretically, the
fast degrees of freedom are driven by the macroscale variables in some concealed
way which we seek to identify through an appropriate model. We are supposed to
have access to at least one measurement set for this modeling task, including
all scales forming the training data set. It is assumed that the resolved scales
control the dynamics through the modulation of the model parameters. We seek to
describe the measured process with a low-dimensional model, recover its temporal
parameter modulation and re-express them by some appropriated characteristic
variable of the resolved dynamical system. Thereby it is assumed that there is a
hidden functional relationship, which we aim to infer, between the model
parameters and the resolved deterministic variable.

The parameter-scaling functions constitute the critical aspect in the identification
of the data-driven stochastic closures. Their existence in real applications
requires exploration, and the functional form is case-specific. Here we test
the concept on synthetically generated examples of different complexity and show
that the scaling functions are deducible.  Our approach provides a way to
develop data-driven stochastic parameterization schemes, focusing on
one-dimensional predictive models.

The paper is organized as follows. Section 2 recalls the theoretical background
on the parameter estimation of a general, stationary SDE and the variational
clustering framework. We explain the coupling of the two inverse problems and
the necessary modifications to the original variational clustering algorithm.
Section 3 covers new approaches in selecting clusters and the regularisation
parameter required to infer the non-stationary models. Section 4 provides a 
summary of our implementation of the clustering method. The numerical examples
in section 5 demonstrate the recovering of the predefined scaling functions, and
section 6 covers the discussion and concluding remarks. \newpage
\section{A concept of data-driven stochastic parameterization}
The concept of data-driven stochastic parameterization is introduced in three
steps. The first two steps form self-contained procedures. The third step is
interweaving the two previous ones and aims at deriving a strategy to develop
data-driven, application-oriented stochastic closure models. The generalized
formulation requires researchers to embed some degree of pre-existing knowledge
about the modeled process in a suitable model structure. After a successful
system identification, a one-dimensional stochastic equation is obtained,
coupled to a macroscale model.

\subsection{Parameter Identification for Autonomous Stochastic Differential 
Equations.}
Consider a stochastic process of interest that
is characterized by a one-dimensional, time-dependent variable $X(t)$. Let 
$X(t_i)$ be the (experimental) observations at observation times $0 = t_0 <
t_1 < \cdots < t_N = T$. For simplicity we assume that the data is recorded at a
constant frequency, so that $t_{i+1} - t_{i} = \Delta t > 0$ for $i= 0, \dots,
N$. The available data set consists of $N$ observations. We consider as underlying data-generating process
a general one-dimensional It$\widehat{\textrm{o}}$ SDE of the form
\begin{equation}\label{eq.:general_sde_form}
	\textrm{d}X = f(X;\bm{\Theta})\textrm{d}t + g(X;\bm{\Theta})\textrm{d}W_t\, ,
\end{equation}
over a finite time interval $[0, T]$, $T > 0$. The functional form for the drift
$f(\cdot;\cdot)$ and diffusion $g(\cdot;\cdot)$ is nonlinear in $X$ and linear
in $\bm{\Theta}$, depending on some unknown vector-valued parameter $\bm{\Theta}
=[\theta_1, \theta_2, \dots ,\theta_n] \in \Omega_{\bm{\Theta}} \subset
\mathbb{R}^n$, where $\Omega_{\bm{\Theta}}$ denotes the admissible parameter
space. The vector $\bm{\Theta}$ contains parameters of both the drift and
diffusion terms. Furthermore, let $W_t = W(t)$ be on admissible, one-dimensional
Wiener process. It is assumed that Eq. \eqref{eq.:general_sde_form} has a unique
strong solution for any $\bm{\Theta}$. The goal is to infer $\bm{\Theta}$ given
$N$ observations using the maximum likelihood approach.  The discrete-time
negative log-likelihood function is defined as \citep{fuchs_inference_2013}
\begin{equation}\label{eq.:nigative_log_likelihood}
	l_N(\bm{\Theta}) := - \sum_{i=0}^{N-1} 
	\textrm{ln}\left[(p_X(\Delta 
	t, X(t_{i+1})|X(t_i); \bm{\Theta}))\right]\, ,
\end{equation}
where  $p_{X}$ denotes the conditional transition density function of the
stationary process $X$ (see Eq. \eqref{eq.:general_sde_form}). From now on we
call $l_N$ as the \textit{fitness function} to match the definition of the
variational framework defined by \cite{horenko_finite_2010}. The optimum of the fitness function with respect to $\bm{\Theta}$ yields the parameter estimate
\begin{equation}\label{eq.:argmin}
\bm{\Theta^{\ast}} \in \underset{\bm{\Theta}\in \Omega_{\bm{\Theta}}}{\textrm{arg min}}\,\,
l_N(\bm{\Theta})\, ,
\end{equation}
which is the MLE. It is known that the MLE converges to the true parameter
$\bm{\Theta}_0$ in probability as $N \rightarrow \infty$ for a fixed $\Delta t$.
Furthermore, it is consistent and asymptotically normal regardless of the
discretization \citep{dacunha-castelle_estimation_1986},
\begin{equation}\label{eq.: asymptotic_normality}
	\sqrt{N}(\bm{\Theta^{\ast}} - \bm{\Theta}_0) \rightarrow^d \mathcal{N}(0, 
	\bm{\sigma}^2_{MLE})\, ,
\end{equation}
where $\bm{\sigma}^2_{MLE}$ is the variance of the MLE
$\bm{\Theta^{\ast}}$.
This estimator converges in distribution $\rightarrow^d$ to the unknown 
parameter at a rate $1/\sqrt{N}$.

In general, the conditional transition density function is not known and to
approximate it we adapt the closed-form expansion following \cite{ait-sahalia_transition_1999, ait-sahalia_maximum_2002}. To obtain the closed-form approximation, the process $X$ is transformed into a process $Z$ in two steps. The transition density of $Z$, namely $p_Z$ is expanded via Hermite polynomials and the needed density $p_X$ is obtained by reverting the transformations.

\begin{figure}[h]
\centering
\includegraphics[width=1.0\linewidth]{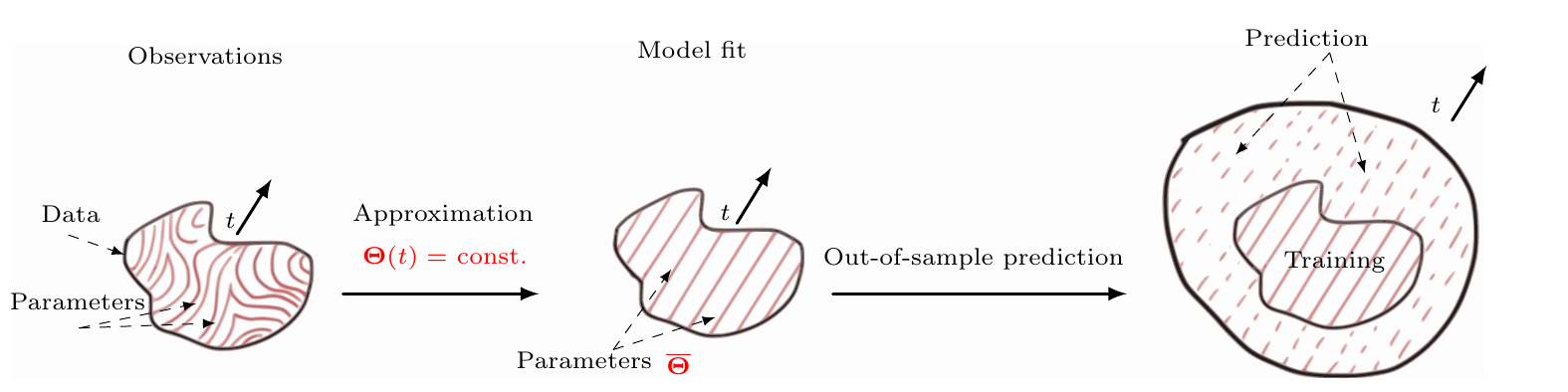}
\caption{An illustration 
of a data-driven modelling concept using SDE with constant parameters.  }
\label{fig:sde_theta_constant}
\end{figure}

The first transformation step $Y:=F(X)$ standardizes the process $X$, such that $Y$ have unity diffusion:
\begin{equation}\label{eq.:transformed_X_in_Y}
	\textrm{d}Y = \mu(Y; \bm{\Theta}) \textrm{d}t + \textrm{d}W_t\, ,
\end{equation}
where the drift $\mu$ is given as
\begin{equation}\label{eq.: Y_drift}
\mu(Y; \bm{\Theta}) = 
\frac{f(F^{-1}(Y);\bm{\Theta})}{g(F^{-1}(Y);\bm{\Theta})} - 
\frac{1}{2}
\frac{\textrm{d} 
	g(F^{-1}(Y);\bm{\Theta})}{\textrm{d} X}\, .
\end{equation}
The Lamperti transform $F$ is
\begin{equation}\label{eq.:lamperty}
F(X) := \int^X g(v, \bm{\Theta})^{-1} \textrm{d}v\, ,
\end{equation}
such that the inverse $F^{-1}(Y)$ exists. The second transformation
is introduced to normalize the increment of $Y$:
\begin{equation}\label{eq.:fix_condition}
	Z:=(Y - y_0) / \sqrt{\Delta t}\, ,
\end{equation}
where $y_0 = Y(t_{i+1}) | Y(t_i)$. The factor $1/\sqrt{\Delta t}$ prevents the 
transition density to become a Dirac delta function for small $\Delta t$ and the 
shift operation sets the initial value for the process $Z$ to $z_0=0$ 
\citep{fuchs_inference_2013}[ch. 6.3]. The truncated Hermite expansion of the transition density function $p_Z$ is defined as
\begin{equation}\label{eq.:expantion}
p_Z^{(J)}(\Delta t, z|y_0; \bm{\Theta}) := \phi(z) \sum_{j=0}^{J} 
\eta_j(y_0, 
\Delta t; \bm{\Theta})H_j(z)\, ,
\end{equation}
where $\phi(z) = e^{-z^2/2}/ \sqrt{2\pi}$ is the weight function, 
\begin{equation}\label{the_coe}
\eta_j(y_0, \Delta t; \bm{\Theta}) := (1/j!) \int_{-\infty}^{+\infty} H_j(z)p_Z(\Delta t, z|y_0; \bm{\Theta})\, dz
\end{equation}
the coefficients and $H_j(z)$ the Hermite polynomials (see Appendix D). By
reverting the transformation from $Z$ to $X$, we obtain the approximated density
$p_X^{(J)}$:
\begin{align}\label{eq.:px transformed}
p_X^{(J)}(\Delta t, X(t_{i+1})|X(t_i); \bm{\Theta})&  = \frac{ 
	1}{g(X(t_{i+1});\bm{\Theta}) 
	\sqrt{\Delta t}}\,  p_Z^{(J)}\left(\Delta t, z |y_0; 
\bm{\Theta}\right)\, .
\end{align}
Theorem 1 of \cite{ait-sahalia_maximum_2002} provides conditions that ensure the convergence of 
\begin{equation}\label{convergence}
p_X^{(J)}(\Delta t, X(t_{i+1})|X(t_i); \bm{\Theta})\rightarrow p_X(\Delta t, X(t_{i+1})|X(t_i); \bm{\Theta}) \qquad \text{as } J \rightarrow \infty\, ,
\end{equation}
both uniform in $X(t_i)$ and $\bm{\Theta}$; resulting in the true transition
density function $p_X$. To evaluate the integral in Eq. \eqref{the_coe}, the
author proposes a Taylor series expansion in $\Delta t$ for the coefficients
$\eta_j$ up to order $M=3$ (see Appendix D). The order of the truncation error
relates to the choice of $J$ as:
\begin{equation}\label{eq.:taylor_truncation_error}
\eta^{(M)}_j(y_0| \Delta t; \bm{\Theta}) = \eta_j(y_0, 
\Delta t; \bm{\Theta}) + \mathcal{O}(\Delta t ^{M+1})\, , \quad j=1, \dots, J;\, \quad M=J/2\, .
\end{equation}
We choose $J = 6$ throughout this work; it is a sufficiently accurate approximation for the proposed clustering approach as has also been advocated by \cite{ait-sahalia_maximum_2002}.
This subsection is closed by providing a summarizing abstract sketch that
represents the scope of this approach where a time series is modeled with an autonomous SDE \eqref{eq.:general_sde_form} (see Fig. \ref{fig:sde_theta_constant}).
In the next subsection, we advance the modelling with a more detailed
approximation.
\subsection{The Non-Parametric Clustering  Methodology}
The previous section covers the parameter estimation of a general autonomous
SDE, assuming the values of its parameters stay constant in time. If we consider
that the data comes from an experiment, then this assumption is more valid for
measurements that are carried out in the laboratory under controlled conditions.
Contrarily, the experiments with uncontrollable conditions are often hard to
reproduce. Consider as an example the atmospheric boundary layer. The duration
of a field measurement campaign may extend over several months, providing a
time series with a wide range of dynamical regimes \citep{ sun_turbulence_2012,
	mahrt_stably_2014, vercauteren_statistical_2019}. It is, therefore, necessary to
cluster the data to be able to isolate and understand the associated dynamics.
If one decides to model the non-stationary process with an SDE, then the
circumstance of time-dependent parameters must be taken into account:
\begin{align}\label{eq.:eq_refine}
	\textrm{d}X = f(X;\bm{\Theta}(t))\textrm{d}t + 
	g(X;\bm{\Theta}(t))\textrm{d}W_t\, \\
	\bm{\Theta} : [0, T] \rightarrow \mathbb{R}^n\, , \qquad \bm{\Theta}(t) = 
	[\theta_1(t),
	\theta_2(t), \dots ,\theta_n(t)] \in \Omega_{\bm{\Theta}}
\end{align}
To estimate the time-varying parameters, we use the non-parametric clustering
approach  with smooth regularization 
introduced by \cite{horenko_finite_2010} (called FEM-H1) and further developed by
\cite{pospisil_scalable_2018}. Here, the term FEM is used to indicate that the finite dimensional projection space used for the regularization is constructed by mimicking the basic principles of the finite element method. Moreover, ''H1'' symbolizes that the finite dimensional projection space is a subspace of $H^1(0,T)$.  Thereby, the task of identifying $\bm{\Theta}(t)$ via FEM-H1 consists of solving two different inverse problems, which are coupled together into one estimation procedure. The first problem consists of identifying fixed parameter values, and the second problem is identifying their location in time. In the following, we focus on the second task. The FEM is adopted to simplify the task complexity in determining the approximation of the time-variable parameter $\bm{\Theta}(t)$. Since the term FEM is commonly associated with partial differential equations, we note that this is not the purpose here.

The problem of parameter estimation (see Eq. \eqref{eq.:argmin}) becomes
ill-posed, if we try to estimate a set of time-dependent parameters, because we try solve a high-dimentional problem.

\begin{figure}[h]
	\centering
	\includegraphics[width=0.5\linewidth]{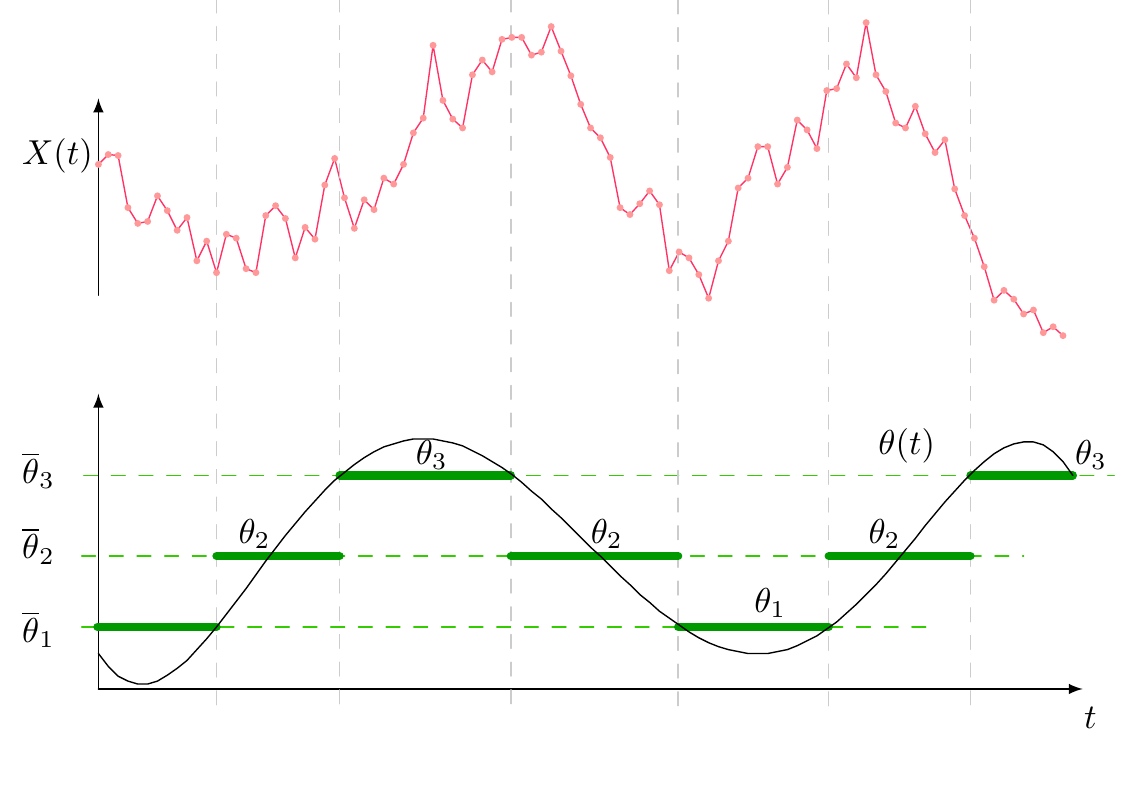}
	\caption{Sketch explaining the approximation for one of the time-dependent 
	parameter functions with 
	piecewise constant subdivisions.}
	\label{fig:approximation}
\end{figure}

To overcome this difficulty, \cite{horenko_finite_2010} suggests regularizing
the fitness function by assuming that the process $X(t)$ varies much faster than
the parameter function $\bm{\Theta}(t)$. Formally, this means expressing the
fitness function with a convex combination of $K$ locally-optimal fitness
functions,
\begin{equation}\label{eq.:horenko_assumption}
	l_N(\bm{\Theta}(t)) := \sum_{i=0}^{N-1} \sum_{k=1}^{K} \gamma_k(t_i)\,  
	f(t_i; \overline{\bm{\theta}}_k)\, , \quad K \geq 2\, ,
\end{equation}
where $f(t_i; \overline{\bm{\theta}}_k)$ is the time-dependent fitness function
formed by the time-averaged parameter value $\overline{\bm{\theta}}_k$. The
overbar denotes the averaged value, which is assumed to be constant. The value
of $K$ indicates the number of clusters or sub-models. The functions
$\gamma_k(t)$ are the affiliation functions with specific properties explained
later. They are unknown and will be estimated with an additional step described
later too. The function $f(\cdot; \cdot)$ is taken to be equivalent to the Eq.
\eqref{eq.:nigative_log_likelihood} and is defined as
\begin{equation}\label{eq.:local_fithess_function}
	f(t_i; \overline{\bm{\theta}}_k) := 
	- \textrm{ln}\left[p_{X}(\Delta 
	t, X(t_{i+1})|X(t_i); \overline{\bm{\theta}}_k)\right]\, , \qquad i = 1, \dots 
	, 
	N-1\, ,
\end{equation}
where the index $i$ iterates over the observations. We keep the sample size of
$f(\cdot; \cdot)$ equal to the sample size of $X$ by setting
the value at the last timestep $N$ to be equal to the first one $f(t_{N};
\overline{\bm{\theta}}_k) = f(t_1; \overline{\bm{\theta}}_k)$. The affiliation
functions $\gamma_k(t)$ contain information on the location of the minima for
each of the $k$-th $f(\cdot; \cdot)$ and are grouped in the so-called
affiliation vector $\bm{\Gamma}(t)=[\gamma_1(t),\gamma_2(t), \dots, \gamma_K(t)]
\in \mathbb{R}^{K}$. The functions $\gamma_k(t)$ satisfy the convexity
constrains forming the partition of unity,
\begin{align}\label{eq.:horenko_cond1}
	\sum_{k=1}^{K}\gamma_k(t) &= 1\, , \qquad \forall t\, \in [0, T],\\
	1 \geq \gamma_k(t) &\geq 0\, , \qquad \forall t \in [0, T], k\, .
	\,\label{eq.:horenko_cond2}
\end{align}
The vector $\bm{\Gamma}(t)$ represents a vector of time-dependent weights. The
functions $\gamma_k(t)$ define the activity of an appropriate $k$-th SDE at a
given time $t$ \citep{pospisil_scalable_2018}, or serves as a
$K$-dimensional probability vector, which expresses a certainty to observe the
$k$-th stationary model at time $t$. For example, assuming that we know the
affiliation vector in advance the temporal evolution of one parameter function
is:
\begin{equation}
	\theta(t) = \sum_{k=1}^{K} \gamma_k(t) \overline{\theta}_k + \varepsilon(t)\, ,
\end{equation}
where $\theta(t)$ is approximated by $K > 2$ piecewise constant values
$\overline{\theta}$ (see Fig. \ref{fig:approximation}) and $\varepsilon(t)$ 
is some approximation error.

As an example, Figure \ref{fig:approximation} illustrates the approximation of
one non-stationary parameter with $K=3$. The figure indicates a variable
location and period of each subdivision which depends on the evolution of
$\theta(t)$ and the number of clusters used. This information is encoded in the
$\bm{\Gamma}(t)$ function. For $K \rightarrow N$ the approximation of
$\theta(t)$ gets better. With each additional cluster used, a different
additional set of parameters need to be determined. If the number of clusters
equals the number of points in the time series ($K=N$), each point should be
potentially described by one equation with its unique parameter values. If, for
example, an equation with $3$ parameters is considered, the time series of
length $N$ is modeled using $3N$ parameters. Such an approximation is considered
to be overfitting. An additional consequence is that the uncertainty for each
estimate of $\overline{\theta}_k$ grows (see Eq. \eqref{eq.:
	asymptotic_normality}) because the amount of data available for each cluster
reduces (i.e. each SDE model).

The averaged clustering functional $L_N$ together with the assumptions
\eqref{eq.:horenko_assumption}, \eqref{eq.:horenko_cond1}, \eqref{eq.:horenko_cond2}
and the fitness function (Eq. \eqref{eq.:local_fithess_function}) becomes:
\begin{equation}\label{eq.:averages_functional}
L_N(\bm{\Gamma}(t), \overline{\bm{\theta}}_1, \dots, 
\overline{\bm{\theta}}_K) = - \sum_{i=0}^{N-1} 
\sum_{k=1}^{K} \gamma_k(t_i)\, \textrm{ln}\left[p_X(\Delta 
t, X(t_{i+1})|X(t_{i}); \overline{\bm{\theta}}_k)\right]\, .
\end{equation}

Applying the regularization and the discretization with the finite element 
method the functional becomes:
\begin{equation}\label{eq.:FEM_discetisation}
L^{\epsilon}_{\widehat{N}}(\bm{\Gamma}(\tau), \overline{\bm{\Theta}}, 
\epsilon^2) 
= 
\sum_{k=1}^{K}\left[ b(\tau, \overline{\bm{\theta}}_k)^{\top} 
\gamma_k(\tau) + \epsilon^2 
\gamma_k(\tau)^{\top} 
\bm{A} 
\gamma_k(\tau) 
\right] \rightarrow \underset{\bm{\Gamma}(\tau), 
	\overline{\bm{\Theta}}}{\textrm{min}}\, .
\end{equation}
where the vector $\bm{\Gamma}(\tau) = [\gamma_1(\tau), \gamma_2(\tau), \dots,
\gamma_K(\tau)] \in \mathbb{R}^{K \times \widehat{N}}$ is the coarsened version of
$\bm{\Gamma}(t)$.	Here, $\tau$ denotes the reduced $t$ grid by a factor $\alpha \in (0, 1]$ and distinguish the notation of the vector $\bm{\Gamma}(t)$ from $\bm{\Gamma}(\tau)$. For additional details on the FEM approximation of $\gamma_k(\tau) $ and the construction of $\bm{A}$ see \cite{horenko_finite_2010} and \cite{metzner_analysis_2012}).

The coarser and uniform mesh with its size $\Delta \tau$ is a set of points $(\tau_j)_{0\leq j \leq \widehat{N}+1}$ with the respective intervals $\mathcal{T}_j = [\tau_j, \tau_{j+1}]$ where $0 = \tau_0 < \tau_1 < \dots < \tau_{\widehat{N}+1} = T$ such that the new number of discretization points $\widehat{N} < N$ . 

The Lagrange $\mathbb{P}_1$ elements associate with the space of globally continuous, locally linear functions for every interval $\Delta \tau$:
\begin{equation}\label{functionspase}
V_{\Delta \tau} = \{v_{\Delta \tau} \in C^0([0,T]): v_{\Delta \tau}|\mathcal{T}_j \in \mathbb{P}_1, 0\leq j \leq \widehat{N} \}\, ,
\end{equation} 
with the subspace:
\begin{equation}\label{subspace}
V_{0,\Delta \tau} = \{ v_{\Delta \tau} \in V_{\Delta \tau}: v_{\Delta \tau}(0) = v_{\Delta \tau}(T) = 0 \}\, .
\end{equation}
The first term
$b(\cdot, \cdot)^{\top} \gamma_k(\cdot)$ in Eq. \eqref{eq.:FEM_discetisation}
is the discretized version of $L_N( \cdot, \cdot)$ (see Eq.
\eqref{eq.:averages_functional}) and is formulated as a dot product. The symbol
$\top$ indicates the transpose of the corresponding vector. The vector
$b(\tau, \overline{\bm{\theta}}_k)$ denotes the discretized model fitness
function,
\begin{equation}\label{eq.:dis_fitness}
b(\tau, \overline{\bm{\theta}}_k) = \left( \int_{\tau_1}^{\tau_2} 
f(t, \overline{\bm{\theta}}_k) v_1(t) 
\textrm{d}\tau, \dots 
, \int_{\tau_{\widehat{N}-1}}^{\tau_{\widehat{N}}} f(t, 
\overline{\bm{\theta}}_k) 
v_{\widehat{N}}(t) 
\textrm{d}\tau \right)\, ,
\end{equation}
where $v(t)$ is the hat function and $f(\cdot, \cdot)$ the function in Eq.
\eqref{eq.:local_fithess_function}. The operation in Eq. \eqref{eq.:dis_fitness}
is denoted as the reduction and it projects a vector into the space of the
finite element. Consequently, the function $f(\cdot, \cdot)$ changes its
sampling rate from $\Delta t$ to a larger mesh size $\Delta \tau$. The second
term $\epsilon^2 \gamma_k(\tau)^{\top} \bm{A} \gamma_k(\tau)$ penalizes the
functional by controlling the regularity of the vector $\bm{\Gamma}(\tau)$,
which is required to obtain a valid solution. With the $\mathbb{P}_1$ elements , the cluster respective stiffness matrix $\bm{H}_k$ assembles to
\begin{equation*}
	\bm{H}_k = 
	\begin{pmatrix}
		2 & -1 & 0 & \cdots & 0 \\
		-1 & 2 & -1 & \cdots & 0 \\
		0 & -1 & 2 & \cdots & 0 \\
		\vdots & \vdots & \vdots & \ddots & \vdots \\
		0 & 0 & 0 & \cdots & 2
	\end{pmatrix} \in \mathbb{R}^{\widehat{N} \times \widehat{N}} \, , \quad k = 1, 
	\dots , K\, .
\end{equation*}
and the total stiffness matrix
\begin{equation*}
\bm{A} = 
\begin{pmatrix}
\bm{H}_1 & 0 & \cdots & 0 \\
0 & \bm{H}_2  & \cdots & 0 \\
\vdots & \vdots & \ddots & \vdots \\
0 & 0 & \cdots & \bm{H}_K
\end{pmatrix} \in \mathbb{R}^{K\widehat{N} \times K\widehat{N}}\, .
\end{equation*}
The regularized variational problem Eq. \eqref{eq.:FEM_discetisation} is solved  
iteratively:
\begin{enumerate}
	\item Fix $\bm{\Gamma}(t)$ and find $\overline{\bm{\Theta}}^{\ast}$ by 
	solving $L_N(\bm{\Gamma}(t), \overline{\bm{\Theta}}) \rightarrow 
	\underset{ \overline{\bm{\Theta}}}{\textrm{min}}$ \eqref{eq.:averages_functional}
	\item Fix $\overline{\bm{\Theta}}$ and find $\bm{\Gamma}^{\ast}(\tau)$ by 
	solving $L^{\epsilon}_{\widehat{N}}(\bm{\Gamma}(\tau), 
	\overline{\bm{\Theta}}, 
	\epsilon^2)  \rightarrow \underset{\bm{\Gamma}(\tau) }{\textrm{min}} $ \eqref{eq.:FEM_discetisation}
\end{enumerate}
The second step is a problem with a quadratic cost function formed by linear
equality and inequality constrains, which is reformulated from Eq.
\eqref{eq.:FEM_discetisation} by defining column vectors from $K$ clusters. The
so-called vector of modeling errors $\bm{B}(\tau)$
\citep{pospisil_scalable_2018} is constructed using the reduced local fitness
function Eq. \eqref{eq.:dis_fitness},
\begin{equation}
	\bm{B}(\tau):=[b(\tau, \overline{\bm{\theta}}_1), b(\tau, 
	\overline{\bm{\theta}}_2 ), \ldots, b(\tau, \overline{\bm{\theta}}_K) ] \in 
	\mathbb{R}^{K\widehat{N}}\, ,
\end{equation}
From Eq. \eqref{eq.:FEM_discetisation} the block-structured 
quadratic programming (QP) problem is:  
\begin{align}\label{eq.:Lucas_problem}
	\bm{\Gamma}^{\ast}(\tau) &:= \textrm{arg } \underset{\bm{\Gamma}\in 
	\Omega_{\Gamma} 
	}{\textrm{min }} 
	L^{\epsilon}_{\widehat{N}}(\bm{\Gamma}(\tau), 
	\overline{\bm{\Theta}}, 
	\epsilon^2)\, , \\
	L^{\epsilon}_{\widehat{N}} &:= \frac{1}{\widehat{N}} \bm{B}^{\top}(\tau) 
	\bm{\Gamma}(\tau) + \frac{1}{\widehat{N}} \epsilon^2 
	\bm{\Gamma}^{\top}(\tau)\bm{A} \bm{\Gamma}(\tau)\, , \\
	\Omega_{\Gamma} &:= \{ \bm{\Gamma}(\tau) \in \mathbb{R}^{K\widehat{N}}: 
	\bm{\Gamma}(\tau) \geq 0\, \wedge\,  \sum_{k=1}^{K} \gamma_k(\tau) = \bm{1} \,  
	\forall t \in [0,T]\}\, ,
\end{align}
where $1/\widehat{N}$ is a scaling coefficient to avoid large numerical values.
To solve the QP, we use the adapted spectral projected gradient method developed
by \cite{pospisil_scalable_2018}. Their method enjoys high granularity of
parallelization and is suitable to run on GPU clusters. Furthermore, it is
efficient and outperforms traditional denoising algorithms in terms of the
signal-to-noise ratio.

Reducing the samples size to $\widehat{N} < N$ reduces the size of the QP problem, and its computational complexity. Higher-order elements also require additional collocation points, but offer smoother approximations. However, the accuracy in determining
$\overline{\bm{\Theta}}^{\ast}$ depends on the number of available observations
as it scales with $1/\sqrt{N}$ and needs, therefore, to be solved at the smallest
available time step $\Delta t$.  If we now reduce the size of the QP problem
(see Eq. \eqref{eq.:dis_fitness}) the accuracy of the estimate
$\overline{\bm{\Theta}}^{\ast}$ will decrease. Consequently, we cannot solve
both problems at the same timestep $\Delta \tau$ and need to interpolate the
solution of QP in each iteration step. This is done as follows:
\begin{enumerate}
	\item Fix $\bm{\Gamma}(t)$ and find $\overline{\bm{\Theta}}^{\ast}$ by 
	solving $L_N(\bm{\Gamma}(t), \overline{\bm{\Theta}}) \rightarrow 
	\underset{ \overline{\bm{\Theta}}}{\textrm{min}}$
	\item reduce the fitness function to the $\tau$-grid 
	\item Fix $\overline{\bm{\Theta}}$ and find $\bm{\Gamma}^{\ast}(\tau)$ by 
	solving $L^{\epsilon}_{\widehat{N}}(\bm{\Gamma}(\tau), 
	\overline{\bm{\Theta}}, 
	\epsilon^2)  \rightarrow \underset{\bm{\Gamma}(\tau) }{\textrm{min}} $ 
	\item interpolate $\bm{\Gamma}^{\ast}(\tau)$ to the $t$-grid step size $\Delta t$
\end{enumerate}

\begin{algorithm}[h]\label{algo.:subspace}
	\caption{The adapted subspace algorithm for SDE models. The main change to the 
	original version of he algorithm is in line 9, where the functional $L_N$ is formulated for a 
	general SDE. Furthermore, lines 14 and 18 are added to maintain the 
	accuracy of the estimated $\overline{\bm{\Theta}}^{\ast}$. The reduced 
	variables are 
	indicated by 
	the time variable $\tau$ and the non-reduced with $t$, respectively. See 
	Appendix A for more details.}
	
	\SetKwInOut{KwIn}{Input}
	\SetKwInOut{KwOut}{Output}
	
	\KwIn{Time series $X$, number of clusters $K$, regularisation value 
	$\epsilon^2$ and the reduction value $\alpha$ }
	\KwOut{$\bm{\Gamma}^{\ast}(t)$ and $ \overline{\bm{\Theta}}^{\ast}$}
	
	Generate random initial $\bm{\Gamma}^0(t)$ satisfying 
	\ref{eq.:horenko_cond1} and \ref{eq.:horenko_cond2}
	
	\tcc{We need the full $\bm{\Gamma}^0(t)$ to find $ \overline{\bm{\Theta}}$ 
	and the reduced $\bm{\Gamma}^0(\tau)$ to solve the QP problem}
	
	Reduce the initial vector $\bm{\Gamma}^0(t)$ to $\bm{\Gamma}^0(\tau)$ 
	on the coarser grid 
	$\Delta \tau$
	
	Set iteration counter $j = 0$ for the main optimization loop
	
	\tcc{main optimization loop}
	\While{$\left| 
	L^{\epsilon}_{\widehat{N}}(\bm{\Gamma}^{j+1}(\tau), 
		\overline{\bm{\Theta}}^{j+1}, 
		\epsilon^2) - 
		L^{\epsilon}_N(\bm{\Gamma}^{j}(\tau), 
		\overline{\bm{\Theta}}^{j}, 
		\epsilon^2) \right|$  $>$ tol}{
	
	\tcc{Estimate $\overline{\bm{\Theta}}^j$ from the full (This will ensure 
	maximum accuracy for the estimate of $\overline{\bm{\Theta}}^j$) 
	$\bm{\Gamma}^j(t)$ 
	applying $K$ times unconstrained minimization }
		\For{$k \leftarrow 1$ \KwTo $K$}{
			$\overline{\bm{\theta}}_k^{j+1} = 
			\underset{\bm{\theta}}{\textrm{arg 
			min}}\, L_N(\gamma_k^j(t), \bm{\theta})$ \tcp*[f]{See eq. 
			\ref{eq.:averages_functional}.}
		}
		\tcc{ From the found $\overline{\bm{\theta}}_k^{j+1}$ compute the 
		fitness functions $f^{j+1}(t, \overline{\bm{\theta}}_k^{j+1})$ or 
		save it from the step in line 9} 
	
		Compute $f^{j+1}(t, \overline{\bm{\theta}}_k^{j+1})$, if not saved 
		from line 9 \tcp*[f]{}
	
		\tcc{Reduce $f^{j+1}(t, \overline{\bm{\theta}}_k^{j+1})$ because 
		$\bm{\Gamma}(\tau)$ will be estimated on the coarser grid 
		$\Delta \tau$}
		
		Compute $b^{j+1}(\tau, \overline{\bm{\theta}}_k^{j+1}) = 
		reduction(f^{j+1}(t, \overline{\bm{\theta}}_k^{j+1}), \Delta 
		\tau, \alpha)$ 
		\tcp*[f]{See eq. \ref{eq.:dis_fitness}}
		
		\tcc{Apply quadratic programming to solve the constrained 
		minimization} 
		
		$\bm{\Gamma}^{j+1}(\tau) = 
		\underset{\bm{\Gamma}}{\textrm{arg 
				min}}\, 
				L^{\epsilon}_{\widehat{N}}(\bm{\Gamma}(\tau), 
				b^{j+1}(\tau, \overline{\bm{\theta}}_k^{j+1}), 
		\epsilon^2) $ \quad satisfying \ref{eq.:horenko_cond1} 
			and \ref{eq.:horenko_cond2}
		
		\tcc{Interpolate $\bm{\Gamma}^{j+1}(\tau)$ onto 
		the 
		original grid with a scale $\Delta t$}
	
		$\bm{\Gamma}^{j+1}(t) = 
		interpolate(\bm{\Gamma}^{j+1}(\tau), \Delta 
		t)$
		
		\tcc{ Addvance the main loop counter j}
		
		$j = j + 1$
	}
	
	\KwRet{$\bm{\Gamma}^{j+1}(t)$, $\overline{\bm{\Theta}}^{j+1}$}
\end{algorithm}\par

In that way, we keep the stepsize of the $\bm{\Theta}$ solver step (1) ($\Delta t$)
separated from the stepsize of the $\bm{\Gamma}$ solver (3) ($\Delta \tau$ ), ensuring
an accurate approximation of $\overline{\bm{\Theta}}$ while keeping the benefit
of the FEM. The algorithm \ref{algo.:subspace} to solve the variational problem
\eqref{eq.:FEM_discetisation} is called the \textit{subspace clustering algorithm}
\citep{horenko_finite_2010} and is reproduced with the proposed modification.

After a successful application of the clustering method, one has several
directions for further research. For instance, the vector
$\bm{\Gamma}^{\ast}(t)$ provides a classifier and can be used to analyze
quantities of interest across the clusters. Instead, we demonstrate a way how
the vector $\bm{\Gamma}^{\ast}(t)$ can be used to construct stochastic closures to effectively model unresolved degrees of freedom as modulated by resolved variables behaviour.

\newpage
\subsection{Stochastic Closure}
One difficulty in finding a stochastic closure is the missing information
regarding the future realization of the vector $\bm{\Gamma}^{\ast}(t)$ because
it is an estimate that is based on the available data. Although we have assumed that we
observe a sufficient amount of  data to capture the full range of regime
dynamics, we are incapable of predicting the future state of the vector
$\bm{\Gamma}(t)$ without additional assumptions and modeling work. One can
associate the affiliation function with a Markov process, assuming that it is
homogeneous and stationary \citep{metzner_analysis_2012}. Alternatively,
\cite{kaiser_stochastic_2017} applied a neural network within the clustering
framework to predict the velocity components of internal gravity waves beyond
the training data set.

We show a different approach by regressing the clustered parameter values to
some known auxiliary process $u(t)$. The limitation is then that it may not
always be possible to find such a process. However, in the case of existence, no
assumptions are required. This approach is attractive since the process $u(t)$
should exist on a larger scale than the process $X(t)$ and thus be predictable.

For simplicity, suppose we are presented with an observation of a   slowly 
evolving, known process $u(t)$ and a fast non-stationary, stochastic process 
$X(t)$ described by some nonlinear SDE. Furthermore, suppose that we have 
access to some expert knowledge, intuition or first principle derivation 
suggesting that $u(t)$ is controlling the non-stationary process $X(t)$ by 
modulating its model parameters. The closure approach is then the following. 
First one estimates the time-varying parameters of the SDE:
\begin{align}\label{eq.:eq_refines}
\textrm{d}X = f(X;\bm{\Theta}(t))\textrm{d}t + 
g(X;\bm{\Theta}(t))\textrm{d}W_t\, \\
\bm{\Theta} : [0, T] \rightarrow \mathbb{R}^n\, , \qquad \bm{\Theta}(t) = 
[\theta_1(t),
\theta_2(t), \dots ,\theta_n(t)] \in \Omega_{\bm{\Theta}}, 
\end{align} 
with the clustering approach presented in section 2.2 that is by solving the 
variational problem:
\begin{equation}\label{eq.:averages_functional_var_prob}
(\bm{\Gamma}^{\ast}, \overline{\bm{\Theta}}^{\ast}) \in 
\underset{\bm{\Gamma} \in \Omega_{\bm{\Gamma}}}{\underset{\bm{\Theta}\in 
		\Omega_{\bm{\Theta}}}{\textrm{arg min}}}\,\, 
		L^{\epsilon}_{\widehat{N}}(\bm{\Gamma}, \overline{\bm{\Theta}}, 
\epsilon^2)\, ,
\end{equation} 
where the number of clusters $K \geq 2$. The functional $L_{\widehat{N}}^{\epsilon}$ is defined in Eq. \eqref{eq.:FEM_discetisation}.

\begin{figure}[h]
	\centering
	\includegraphics[width=1.0\linewidth]{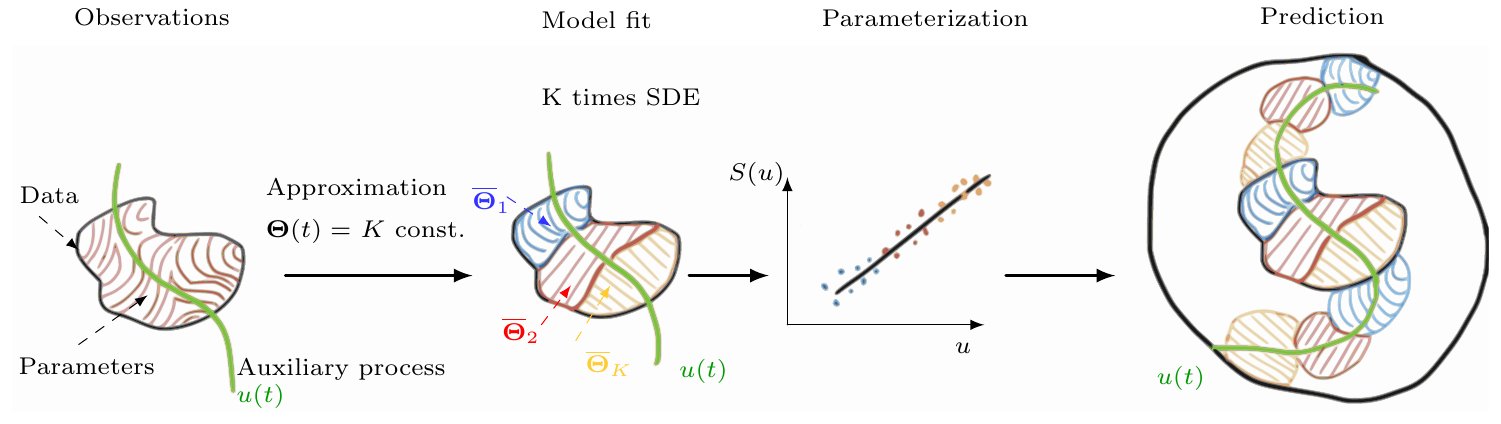}
	\caption{The modeling approach with FEM-H1 regularization of SDE. The 
		temporal modulation of the parameters is parameterized through the slow 
		scale 
		process $u$. The functional relation $S(u)$ is hidden and needs to be 
		discovered a posteriori. }
	\label{fig:step3}
\end{figure}

Second, with the help of $\bm{\Gamma}^{\ast}(t)$ one estimates the cluster respective values of the process $u(t)$ by computing the weighted average:
\begin{equation}\label{eq.:weighted_arithmetic_mean}
\overline{u}_k = \frac{\sum_{i=1}^{N}  
	u(t_i) \gamma_k(t_i)}{\sum_{i=1}^{N} \gamma_k(t_i)}\, , \qquad k = 1, \dots 
, K \, .
\end{equation}
where $\overline{u}_k$ denotes the averaged value of the auxiliary process
$u(t)$ in cluster $k$. Recall that the $\gamma_k(t_i)$ is the affiliation
function for cluster $k$ (see Fig. \ref{fig:OU_example_result}a the color coded
regions). In the last step, one performs a regression analysis to parameterize
the parameter variability with some appropriate functions $S_m$, which map the
values of $u(t)$ to $\theta_{m}(t)$:
\begin{equation}
	S_m: \overline{u}_k \rightarrow \overline{\bm{\theta}}_{k, m}\, , \qquad 
	\overline{\bm{\theta}}_k = [\overline{\theta}_{k,1}, 
	\overline{\theta}_{k,2}, \dots, \overline{\theta}_{k,n}], \qquad k = 1, 
	\dots , K \, \qquad\, m = 1, \dots , n .
\end{equation}
There is no restriction on the differentiability of the functions $S_n$, it can 
be
parameterized by a piecewise-defined function. The parameterization is performed
for each element of the parameter vector $\bm{\Theta}$ separately. The vector of
scaling functions $\bm{S} = [S_1, S_2, \dots, S_n] \in \mathbb{R}^n$ defines the
modulation of the model parameters through some process $u(t)$. The complexity
of the relation is application dependent and needs further exploration. If all
scaling functions can be found, we obtained a closed system to predict the
process $X(t)$
\begin{equation}\label{eq.:closed_system}
\textrm{d}X = f(X;\bm{S}(u(t)))\textrm{d}t + 
g(X;\bm{S}(u(t)))\textrm{d}W_t\, ,
\end{equation}
where $u(t)$ is known for every $t$. The approach allows to model
non-equilibrium processes which are characterized by the nonlinearity in $f$
and $g$ in the cases where the nonstationarity of the process is difficult to
understand. Figure \ref{fig:step3} summarizes the idea of this concept. In
practice the functional relation between process $u(t)$ and the modulated
parameters of the SDE $\bm{\Theta}(t)$ is assumed to hold for the
estimated parameter range, as shown in Section 5.
\pagebreak

 \clearpage
\newpage
\section{Selection of the Hyperparameters}
The selection  of a suitable model includes the determination of a suitable
regularization parameter $\epsilon^2$ , controlling the smoothness of the
affiliation vector, the number of clusters $K$ and the functional forms
$f(\cdot; \cdot)$ and $g(\cdot; \cdot)$.  A variety of strategies with examples
are described in \citep{horenko_finite_2010, horenko_identification_2010,
	metzner_analysis_2012} for example, which are partially based on identifying
several criteria simultaneously. Here we present new approaches to select
$\epsilon^2$ and $K$ independently from each other, providing a robust
estimation strategy. In particular our approach for selecting the parameter
$\epsilon^2$ accounts indirectly for the frequency content of the auxiliary
process $u(t)$.

\subsection{Optimal Regularization Parameter}\label{sec.:opt_regularisation}
\begin{figure}[h]
\centering
\includegraphics[width=0.7\linewidth]{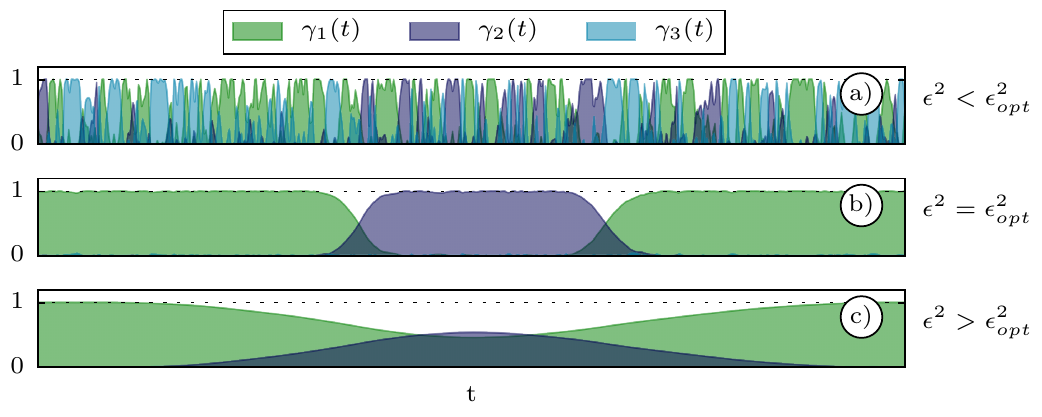}
\caption{An example on the different regularity of the affiliation vector 
$\bm{\Gamma}^{\ast}(t)$. Each panel displays 
the same time window but for different 
values of $\epsilon^2$. a) The weak regularization: the regime 
periods are irregular. Their temporal scale is close to  $\textrm{d}t$. b) 
The optimal regularization: $\gamma_1(t)$ and 
$\gamma_2(t)$ shows relatively sharp and clear separation. c)The strong 
regularization: $\gamma_2(t)$ is blurred and suppressed.}
\label{fig:gamma_regularity}
\end{figure}

\cite{horenko_finite_2010}  analyzed the impact of the parameter $\epsilon^2$ on
the regularity of the cluster affiliation function $\bm{\Gamma}^{\ast}(t)$. His
findings show that the optimal parameter $\epsilon^2_{\text{opt}}$ is
characterized by a sharp separation between each of the functions $\gamma_k(t)$,
such that the regime transitions are as short as possible. Contrary, poor
partitioning is present if two or more functions $\gamma_k(t)$ occupy the same
period. We recall, that the functions $\gamma_k(t)$ obey the convexity condition
(see Eq. \eqref{eq.:horenko_cond1}; \eqref{eq.:horenko_cond1}) and if several
are active simultaneously then none of the model parameters fits the data
sufficiently well. Such results indicate that the number of clusters is set too
high or the model structure is not suitable for the data under consideration. A
key aspect of choosing a proper value for $\epsilon^2$ is in monitoring the
regularity of the affiliation function $\bm{\Gamma}^{\ast}(t)$. Figure
\ref{fig:gamma_regularity} shows a portion of $\bm{\Gamma}^{\ast}(t)$ with $K=3$
for a toy example.

\begin{figure}
\centering
\includegraphics[width=0.7\linewidth]{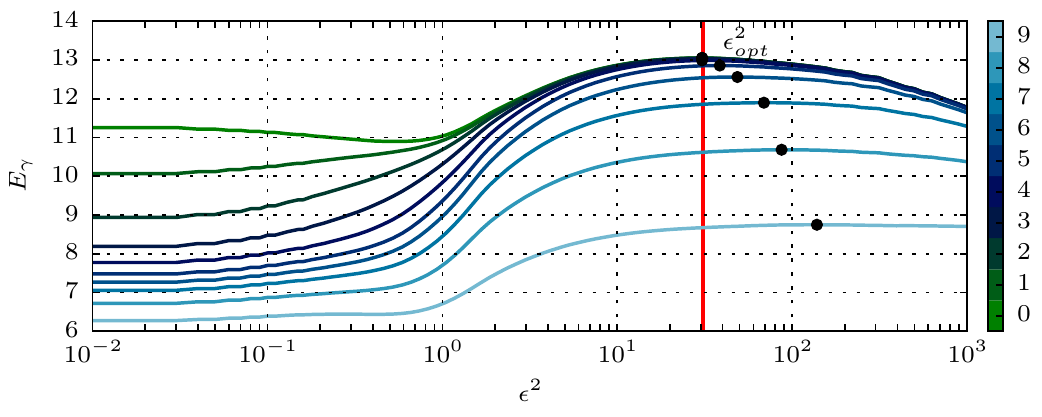}
\caption{An approximation of the optimal regularisation parameter $\epsilon^2$ 
for 
a toy example with three clusters. 
The y-axis is 
the cluster averaged signal energy of $\bm{\Gamma^{\ast}}(t)$ (see Eq. 
\ref{eq.:gamma_averaged_energy_wl_max_dis}). By wavelet filtering the estimated 
affiliation vector we emphasize the maximum. 
The colorbar labels the 
	number of removed details 
	coefficients levels ($9$ of maximum $16$ is shown). The zero marks the 
	unfiltered $\bm{\Gamma^{\ast}}(t)$. In this 
	example one can exclude four levels to emphasise the maximum without 
	compromising the optimum.}
\label{fig:gamma_opt}
\end{figure}

To motivate the strategy in finding $\epsilon^2_{opt}$ we consider
$\bm{\Gamma}^{\ast}(t)$ for different values of the parameter (see Fig.
\ref{fig:gamma_regularity}a, b, c). Note the white regions between $0$ and $1$
in Fig. which are not color coded. The total size of this area is changing in
each of the panels and is approximately the smallest
at the value of $\epsilon^2_{opt}$. If we consider the function
$\bm{\Gamma}^{\ast}(t)$ as an oscillating signal, then to minimize the white
area we search for a $\bm{\Gamma}^{\ast}(t)$ that has the largest energy, where
energy is understood as $||\gamma_k||^2_{L^2(0,T)}$. Using the cluster averaged 
signal energy $E_{\gamma}$ for a given value $\epsilon^2$ 
\par
\begin{equation}\label{eq.:gamma_averaged_energy}
	E_{\gamma}(\epsilon^2) = \frac{1}{K} \sum_{k=1}^{K} \int_{0}^{T} 
	\gamma_{k}^2(t, \epsilon^2) \textrm{d}t\, ,
\end{equation}
where $\gamma_{k}^2(t, \epsilon^2)$ is the solution of Eq.
\eqref{eq.:FEM_discetisation}, we determine the value $\epsilon^2$ that maximize
the energy. The maximum of the curve $E_{\gamma}(\epsilon^2)$ indicates the
optimal regularization value (see Fig. \ref{fig:gamma_opt} color code value
$0$). The true affiliation function which was used to generate the training data
is generated with sharp transitions and three clusters.  Furthermore, in the
identification step, we used the same functional form that was used in the
construction process. In other words, the considered examples are well-posed in
the sense that the data-generating model is perfectly identifiable. However,
working with real, unexplored data we, unfortunately, may not know the true
functional form of the SDE. This may lead to a situation where the functions
$\gamma_k(t)$ are uncertain which affects the curve $E_{\gamma}(\epsilon^2)$ by making it less concave and hence loosing a pronounced maximum. This is an indication of an
increased sensitivity. To achieve a more robust procedure we investigate the
frequency content of the vector $\bm{\Gamma}^{\ast}(t)$.

One way to emphasize the maximum of the curve $E_{\gamma}(\epsilon^2)$ is by
gradually removing the high frequency oscillations in $\bm{\Gamma}^{\ast}(t)$.
Since the affiliation function is not periodic we filter it with a wavelet method \citep{foufoula-georgiou_wavelets_1994}. The affiliation vector is decomposed using the discrete
wavelet transform to obtain $C$ levels with $C+1$ frequency bands. The function
$\gamma_k(t)$ is then represented as  the sum of the low-frequency component
$\gamma_{k, C}^A(t)$ after $C$ levels of transformations and  the sum of the
high-frequency components $\gamma_{k, c}^D(t) $ over the previous decomposition
levels
\begin{equation}\label{eq.:decomposition}
	\gamma_k(t) = \gamma_{k, C}^A(t) + \sum_{c = 1}^{C} \gamma_{k, c}^D(t)\, .
\end{equation}
As a basis for the filtering, we choose the Haar wavelet. They represent the
sharp jumps of the affiliation vector in the best way, especially when
approaching the value $\epsilon^2_{opt}$. Utilizing the $L_2$-orthogonality
property of the Haar wavelets one inserts Eq. \eqref{eq.:decomposition} in to
\eqref{eq.:gamma_averaged_energy} to obtain
\begin{equation}\label{eq.:gamma_averaged_energy_wl}
E_{\gamma}(\epsilon^2) = \frac{1}{K} \sum_{k=1}^{K} \int_{0}^{T} 
\left[ \gamma_{k, C}^A(t)\right]^2 \textrm{d}t +  \frac{1}{K} \sum_{k=1}^{K} 
\int_{0}^{T} 
\left[ \sum_{c = 1}^{C} \gamma_{k, c}^D(t) \right]^2 \textrm{d}t \, .
\end{equation}
The signal $\gamma_{k}(t)$ is decomposed with the maximum number of possible
levels $C$, which is dependent on the number of samples. The first term in Eq.
\eqref{eq.:gamma_averaged_energy_wl} is the mean value of the total time series
and is irrelevant since we are interested in regime transitions between
different models. In the following steps we keep the same notation for
$E_{\gamma}(\epsilon^2)$. The local-support property of the wavelets allows to
take out the summation sign of the integral,
\begin{equation}\label{eq.:gamma_averaged_energy_wl_max}
E_{\gamma}(\epsilon^2) =  \frac{1}{K} \sum_{k=1}^{K} 
\sum_{c = 1}^{C} \int_{-\infty}^{\infty} 
\left[ \gamma_{k, c}^D(t) \right]^2 \textrm{d}t \, .
\end{equation}
The sample size of the discrete-time series $\gamma_{k, c}^D(t)$ is extended to
the nearest value $2^i$ with the reflect-padding method, where $ i \in
\mathbb{N^+}$. Next, we rewrite the continuous integral as a sum
\begin{equation}\label{eq.:gamma_averaged_energy_wl_max_dis}
E_{\gamma}(\epsilon^2) =  \frac{1}{K} \sum_{k=1}^{K} 
\sum_{c = 1}^{C} \sum_{d=1}^{C_d}
\left[ \tilde{\gamma}_{k, c}^D(d) \right]^2 \Delta t \, ,
\end{equation}
where $C_d(c)$ denotes the number of coefficients at the decomposition level 
$c$.
Due to the efficient implementation of the multi-resolution decomposition the
number of coefficients at each level is consecutively reduced by $2$. With each
decomposition level the cut-off frequency is halved and the energy of the
resulting frequency bands can be removed by setting the corresponding
coefficients to zero.

The optimal cut-off frequency may be estimated in several ways. A reasonable
choice would be to set it to the scale of a spectral gap, if one is present,
between the process $X(t)$ and $u(t)$. A spectral gap is formed when the highest
frequency of $u(t)$ is larger than the lowest frequency of $X(t)$. In this way,
the scale of the process $u(t)$ provides a limit for the smallest scale of the
regime jumps in $\bm{\Theta}(t)$. There is no need to define the upper bound
because the curve $E_{\gamma}(\epsilon^2)$ drops naturally due to the
over-smoothing of the vector $\gamma_k(t)$ at higher values $\epsilon^2$. The
over-smoothing is not detected by the sharp jumps of the Haar wavelet basis.
Consequently, the wavelets stop registering the energy at the respective scales
(see the roll-off in Fig. \ref{fig:gamma_opt} for $\epsilon^2 > 10^2$).
Alternatively, one plots the energy curves by consecutively removing the high
frequencies and looks for a magnification of a maximum without a change in its
position (see Fig. \ref{fig:gamma_opt}). The multi-resolution analysis is
performed using the python library PyWavelets \citep{lee_pywavelets_2019}.

\subsection{Optimal Number of Clusters K}\label{sec.:opt_number_of_clusters}
The approach to determine the optimal number of clusters $K_{\text{opt}}$ is
borrowing ideas from the elbow method which is used in the K-means clustering
approach \citep{yuan_research_2019}. The idea is to monitor an error measure for
the clustering procedure over the number of clusters and search for a
characteristic ''kink'' in the graph. The choice of the error measure is,
therefore, important and specific to the method used. For instance, a measure of
compactness in K-means clustering is defined as the sum of the Euclidean
intra-cluster distances between points in a given cluster. We introduce a
different error measure because the SDE clustering method differentiates the
temporal dynamics within the clusters. Two clusters may strongly overlap but be
distinct in terms of their induced stationary probability density functions. The
Euclidean distance would be a poor choice in that case. Instead, we exploit the
Kullback-Leibler (KL) divergence, which is closely related to the probability
density of the SDE:
\begin{equation}\label{eq.:KL_div}
	KL(P||Q) := \int P(x)\,  \textrm{log}\left( \frac{P(x)}{Q(x)} \right) 
	\textrm{d}x\, .
\end{equation}
The KL divergence describes the information gain by changing from the
distribution with density $P$ to the distribution with density $Q$. Here we will
use the $P$ and $Q$ as the stationary distributions of two different clusters.
We seek to measure the difference between $P$ and $Q$ or more generally the
\textit{diversity} between $K$ stationary distributions. Due to the asymmetry of
KL divergence we consider also the values of the opposite divergences
$KL(Q||P)$. Let $p^{\star}_k$ denote the stationary distribution of an SDE with
the estimated parameter set $\overline{ \bm{\theta}}^{\ast}_k$ and the
corresponding $\gamma_k(t)$  for cluster $k = 1, \dots, K$. For each cluster we 
thus have a different SDE model
that provides a stationary distribution. We define the \textit{weighted
	diversity matrix} $D = (d_{ij}) \in \mathbb{R}^{K \times K}$ computed as
\begin{equation}\label{eq.:dist_matrix}
	d_{ij} := KL(p_i^{\star}||p_j^{\star})= \nu_j \int 
	p_i^{\star}\,  
	\textrm{log}\left( \frac{p_i^{\star}}{p_j^{\star}} \right) 
	\textrm{d}x\, \geq 0 \qquad i,j = 1, \ldots \, K\, ,
\end{equation}
where $\nu_j$ are the weights computed from the  affiliation vector that 
quantify the relative frequency of cluster $j$:
\begin{equation}\label{eq.:weights}
\nu_j= \int_{0}^{T} \gamma_j(t)\textrm{d}t\, , \qquad j = 
1, 
\ldots K\, .
\end{equation} 
The weights are introduced for the following reason. By adding more clusters the
less probable ${\gamma_j}(t)$ is suppressed (see for example Fig. \ref{fig:gamma_regularity}c) and barely reaches values equal to
$1$ or it has a short lifetime relative to the full-time horizon $T$. These
suppressed, less-probable, inactive regimes should therefore include less
diversity and compensate the measure. Since we are looking to find the optimal 
number of cluster within some reasonable range $2 < K < K_{\text{max}}$ we 
define our diversity measure of clustering with $K$ clusters as:
\begin{equation}\label{eq.:diversity_measure}
	W_K= \sum_{i = 1}^{K} \sum_{j = 1}^{K} d_{ij}\, , \qquad K = 2,3, 
	\ldots, K_{\textrm{max}}\, ,
\end{equation}
where $K_{\textrm{max}}$ denotes the maximum considered clusters, assuming that
$2 < K_{\textrm{opt}} < K_{\textrm{max}}$. We reiterate that $p^{\star}_k$ in Eq. \eqref{eq.:dist_matrix} denotes the stationary distribution of the
SDE model characterized by parameter set $\overline{ \bm{\theta}}^{\ast}_k$. 
Specifically, for the SDE model 
\begin{equation}\label{eq.:eq_refines2}
	\textrm{d}X = f(X;\overline{ \bm{\theta}}^{\ast}_k))\textrm{d}t + 
	g(X;\overline{ \bm{\theta}}^{\ast}_k))\textrm{d}W_t\, , \qquad k = 1, 
	\dots, K\, ,
\end{equation} 
the invariant density $p^{\star}_k$ can be written as
\begin{equation}\label{eq.:}
p^{\star}_k(X)= \frac{N_c}{g^2(X; \overline{ \bm{\theta}}^{\ast}_k)} 
\textrm{exp}\left(\int 
\frac{2 f(X; \overline{ \bm{\theta}}^{\ast}_k)}{g^2(X; \overline{ 
\bm{\theta}}^{\ast}_k)}dx\right)\, , \qquad k = 1, \ldots K\, ,
\end{equation}
where the normalisation constant $N_c$ is:
\begin{equation}\label{eq.:stationary_pdf_normalisation_constant}
N_c^{-1}= \int \frac{1}{g^2(X; \overline{ \bm{\theta}}^{\ast}_k)} 
\textrm{exp}\left(\int 
\frac{2 
	f(X; \overline{ \bm{\theta}}^{\ast}_k)}{g^2(X; \overline{ 
	\bm{\theta}}^{\ast}_k)}dx\right)\, , \qquad k = 1, \ldots K\, ,
\end{equation}
provided all terms are well-defined (see \cite{horsthemke_noise_1984}[ch. 6.1] 
for details).

In real-world applications the gain in the diversity $W_K$ can be too
steady making the detection of $K_{\textrm{opt}}$ difficult. We use further
extension of the elbow method idea to make the estimation of $K_{\textrm{opt}}$
robust. Indeed, we follow \cite{tibshirani_estimating_2001} where it has been
suggested to measure the increase of $W_K$ with respect to some null reference
distribution. Namely, a data set which has no obvious clustering. In the work of
\cite{tibshirani_estimating_2001} the randomly spreading of points in
2-dimensional space serves as a null reference data set, which is suitable for
clustering with the K-means method. By considering the properties of our
clustering method we choose the Wiener process with reflecting boundaries as the
null reference time series. The optimal number of clusters is estimated as the
point where we get the smallest distance between the clustering-diversity
measure of the analyzed and the reference data:
\begin{equation}\label{eq.:gap}
	\textrm{Gap}(K)= \mathbb{E}[\textrm{log}(W_K^{\ast})] -  
	\textrm{log}(W_K)\, ,
\end{equation}
where $W_K^{\ast}$ denotes the diversity measure of the null reference dataset
and the extracted value $ \mathbb{E}[\textrm{log}(W_k^{\ast})]$ is estimated from clustering $B \gg 1$ different, independent reference time series. The
boundaries for the reference Wiener process is confined to the range of the analyzed process. The section 5 demonstrates the application of presented approach on controlled numerical examples.

 \section{Implementation}
We have developed our version of the algorithm \ref{algo.:subspace} in the
programming language C++. The core components of the FEM-H1 framework, namely: the
spectral projected gradient method for the quadratic programming, the FEM
reduction and interpolation procedures, and the procedure to generate an initial
$\bm{\Gamma}(t)$ are reimplemented from the MATLAB code that is created by original authors
\cite{pospisil_scalable_2018}. For the solution of the $\bm{\Theta}$ problem (see
algorithm \ref{algo.:subspace} line 9) we use the NLopt nonlinear-optimization
package  \citep{g_johnson_nlopt_2021}. It is an interface of many global and
local optimization algorithms which can be tested just by changing one
parameter. Specifically, we applied a combination of Controlled Random Search
(CRS) with local mutation \citep{price_controlled_1977, price_global_1983,
	kaelo_variants_2006} and the "PRAXIS" gradient-free local optimization via the
"principal-axis method" \citep{richard_p_algorithms_2013}. Our code is GPU
accelerated using the CUDA libraries cuSPARSE, Thrust and requires the CUDA
version 10.2. The transition density is computed with the MATLAB symbolical
toolbox and exported to a C++ code. The template library Thrust made it possible
to include the machine-generated code of the transition density $p_X$ into the main program flow. \newpage
\section{Numerical examples}
The numerical study consists of three examples with varying complexity. The
auxiliary process $u(t)$ here is stochastic with a smooth sample path and a
reasonable scale separation between $u(t)$ and $X(t)$. However, in real-world
applications the process $u(t)$ represents some characteristic quantity of a
dynamical system, with the assumption to control the time-varying parameters of
the process $X(t)$.

The numerical examples consist of the step by step construction of a synthetic
time series $X(t)$ and the identification of the underlying functional
relationship between the parameters of a predefined SDE and some random process
$u(t)$. Here it is assumed that the process  $u(t)$ is known at any instance of
time and it serves as a predictor for the time-varying parameter vector
$\bm{\Theta}(t)$. To create the synthetic data set $X(t)$, the process $u(t)$ is
used together with the arbitrary created functions $S_n(u)$ to generate the
evolution of the parameter vector $\bm{\Theta}(t) = [S_1(u(t)), S_2(u(t)),
\dots, S_n(u(t))]$. The sample path $X(t)$ is simulated numerically with the
predefined SDE, where the parameters are given by $\bm{\Theta}(t)$. The
time series $X(t)$ forms the data set of the numerical tests, from which we
recover the functions $S_n(u)$ and compare them to the one that have been used in
the creation step.

\subsection{The Underlying Auxiliary Process Model}
As a way of defining an auxiliary process that respects the required scale separation between $X(t)$ and $\bm{\Theta}(t)$, we provide a brief explanation of a stochastic dynamical system to simulate a slow-evolving underlying process. A realization of the auxiliary process $u(t)$ is constructed for all
examples following the same principle. We utilize the 4-dimensional Ornstein -
Uhlenbeck process (OU) supplementing it with the coefficients of the normalized
Butterworth polynomials,
\begin{equation}\label{eq:OU_U(t)}
\textrm{d}\bm{U} = - \frac{1}{T_c} \bm{A}_D \bm{U} \textrm{d}t +  \bm{B} 
\textrm{d}\bm{W}\, , \qquad \bm{U}=[u(t), u_1(t), u_2(t), u_3(t)] \, ,
\end{equation}
where $T_c > 0$ controls the cut off frequency, $\bm{B} > 0$ is the 
diffusion
coefficient. The auxiliary process $u(t)$ is then the first component of the vector $\bm{U}(t)$ The diffusion matrix $\bm{B}$ is almost everywhere zero, and
together with the drift matrix $\bm{A}_D$ takes the form:
\begin{equation}
\bm{A}_D = 
\begin{bmatrix}
	0 & 1 & 0 & 0\\
	0 & 0 & 1 & 0\\
	0 & 0 & 0 & 1\\
	a_0 & a_1 & a_2 & a_3
\end{bmatrix}\, , \qquad
\bm{B} = 
\begin{bmatrix}
0 & 0 & 0 & 0\\
0 & 0 & 0 & 0\\
0 & 0 & 0 & 0\\
0 & 0 & 0 & b_0
\end{bmatrix}\, .
\end{equation}
The coefficients in the matrix $\bm{A}_D$ are: $[a_0, a_1, a_2, a_3] = [1, 2.61,
3.41, 2.61]; b_0 = 1$. The effect of the coefficients in $\bm{A}_D$ is the
smoothening of the first component in the vector $\bm{U}$, which is an analogy
to a filtering approach.	The matrix $\bm{A}_D$ is constructed based on the filter design theory
utilizing the coefficients of the normalized Butterworth polynomials. The
drift term in Eq. 5.1 acts as a filter on the diffusion process
that operates as an input signal. Consider that the process $u(t)$ can also be
created differently. One solves a one-dimensional OU process and then
filters the result with a linear filter by defining the cut-off frequency
$T_c$. Due to the filtering, we obtain a sufficiently smooth time series. The
same result is obtained by solving Eq. 5.1. Alternatively,  one
could prescribe a power spectrum density in the frequency domain and then
transform it into the physical space with the inverse Fourier transform. The
approach above is simple and favorite because we solve for $u(t)$ and $X(t)$
simultaneously and have control of the timestep. However, an analytical
approach would be the sophisticated solution.
\clearpage
\subsection{The Non-Stationary Ornstein - Uhlenbeck Process} 
\label{sec.:ou_example}
We start the study with a trivial case where the sample path $X(t)$ is modeled
by an OU process, which has a linear drift term and additive noise. Consider the
OU process of the following form:
\begin{equation}\label{eq:OU}
\textrm{d}X = (\theta_1(t) - \theta_2(t) X ) \textrm{d}t +  \theta_3(t) 
\textrm{d}W\, , \qquad X(t_0) = 0\, ,
\end{equation}
where $\bm{\Theta}(t) = [\theta_1(t), \theta_2(t), \theta_3(t)]$ will be
modulated by  $u(t)$ defined with Eq. \eqref{eq:OU_U(t)}. The analytical
expression for the transition density $p_X(\cdot, \cdot)$, which is required to
compute the MLE, in this case, is known. However, we still use the approximation
with the Hermite expansion to test the implementation of the framework.
\begin{figure}[h]
\centering
\includegraphics[width=1.0\linewidth]{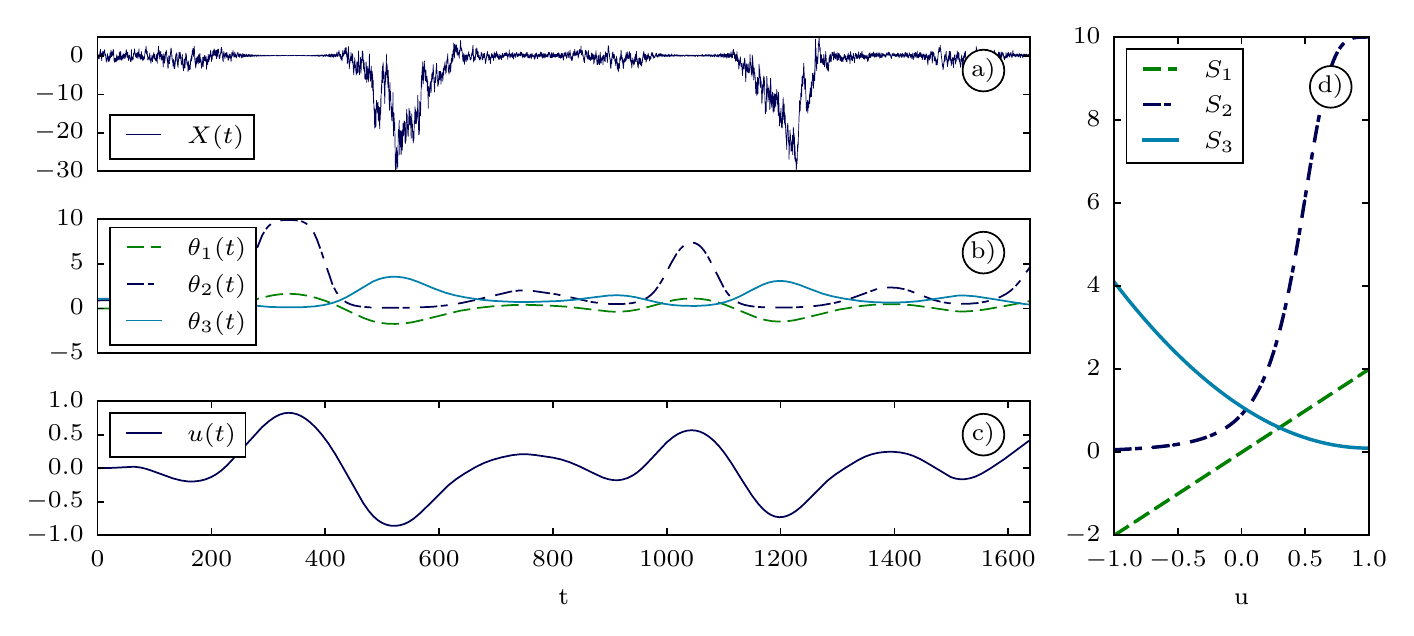}
\caption{The summary of functions used to generate the non-stationary training data 
$X(t)$
according to the Eq. \ref{eq:OU}. a) The sample path of the process 
$X(t)$. b) The temporal evolution of the model parameters. c) The Parameter 
auxiliary process 
$u(t)$. d) The scaling functions $\theta_n(t) = S_n(u(t))$}.
\label{fig:OU_example}
\end{figure}
Figure \ref{fig:OU_example}a demonstrates one realization of the sample path for
the process $X(t)$. The SDE \eqref{eq:OU} is solved simultaneously with the Eq.
\eqref{eq:OU_U(t)} using different independent realizations of the Wiener
process for each of the equations. We use the Euler-Maruyama method with a
time-marching step size of $\Delta t = 10^{-6}$ and then downsample the results
to the time step $\Delta t = 10^{-1}$, which then has in total $16384$ sample
points. The scaling functions are defined as (see fig. \ref{fig:OU_example}d)
\begin{align}
	\theta_1(t) := S_1(u)&= 2 u(t)\, ,\\
	\theta_2(t) := S_2(u)&= \left[(u(t)-1)^4 + 0.1\right]^{-1}\, , \\
	\theta_3(t) := S_3(u)&= (u(t)-1)^2 + 0.1\, .
\end{align}

The evolution of the vector $\bm{\Theta}(t)$ (see fig. \ref{fig:OU_example}b)
and its functional relationship to the process $u(t)$ is hidden and supposed to
be unknown in real applications. In this example, one can observe an
interrelationship between the variance of the process $X(t)$ and the value of
the process $u(t)$ by only considering the sample paths because the equation is
linear. As shown in the next examples, this is not obvious with a nonlinear
system. In the following, we demonstrate how the methodology introduced here is
applied to recover the scaling functions $S_n$, given only the discrete-time
values of $X(t)$ and $u(t)$. For details on running the optimization algorithm
for this example see Appendix B.

\begin{figure}[h]
	\centering
	\includegraphics[width=1.0\linewidth]{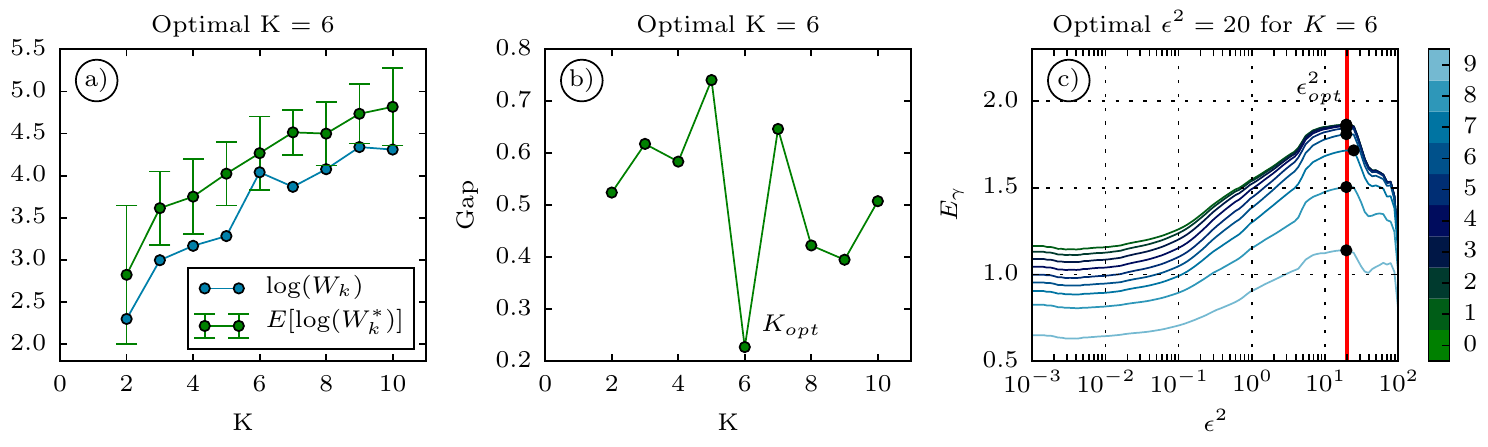}
	\caption{An estimation of the number of clusters for the OU example (see 
	Fig. \ref{fig:OU_example} and Eq. \eqref{eq:OU}) using the adapted gap statistics 
	approach see Sect 3. In 
	panel a) we find that for $K>6$ the diversity of the training data is 
	increasing at the comparable rate as in the reference data sets. This 
	suggest 
	the 
	optimal $K$ to be 6. This is confirmed  by the minimum value in panel b).   
	c) The estimation of the regularization parameter $\epsilon^2$ for the 
	clustering 
	with $K=6$. The colorbar codes the consecutive suppression of high 
	frequency bands to emphasise the  maximum. The effect is minimal.}
	\label{fig:ougammaoptkopt}
\end{figure}
By plotting the cluster averaged signal energy of the affiliation vector, i.e
its  $L_2$-norm, against the parameter $\epsilon^2$ we estimate the optimal
value $\epsilon^2_{\textrm{opt}}$ for each value of $K$ individually (see one
example for $K=6$ in fig. \ref{fig:ougammaoptkopt}c). According to the figure
\ref{fig:ougammaoptkopt}c the maximum energy is found at a value
$\epsilon^2_{\textrm{opt}} \approx 20$, which indicates the clearest separation
of the regimes. As demonstrated in section \ref{sec.:opt_regularisation} the
maximum in the curve $E_{\gamma}$ may be further emphasized by suppressing the
high frequencies in the vector $\bm{\Gamma}^{\ast}(t)$. This is highlighted in
figure \ref{fig:ougammaoptkopt}c with different colors, where index $0$ is the
unfiltered result and each next integer marks the consecutive reduction of the
frequency content by a factor of two (discrete wavelet transform). We find that
the frequency suppression shows the desired effect, although this is unnecessary
in this example because the maximum is well pronounced. Additionally, we have
inspected  the  affiliation function for every value of $\epsilon^2$ (not shown)
and found relatively small differences in its evolution for a wide range of the
value $\epsilon^2$. This indicates a robust estimation and may be attributed to
the simplicity of the example.

The procedure to estimate the optimal $K_{\textrm{opt}}$ requires more effort
than it is  the case with $\epsilon^2_{\textrm{opt}}$. We introduced the
diversity measure in section \ref{sec.:opt_number_of_clusters} to characterize
the incremental information gain with every additional cluster. Figure
\ref{fig:ougammaoptkopt}a demonstrates the graph of this measure (blue) for the
training data $X(t)$ and for the reference data (green). The reference process
is generated $B = 20$ times with the same timestep that was used in the
identification step. The clustering is repeated on each of these data sets with
the value $\epsilon^2_{\textrm{opt}} \approx 20$ for each of the $K$ clusters.
The intention is to keep all parameters of the framework unchanged and only to
replace the data set. The resulting graph of the information gain for the
clustering of  the reference time series is shown in figure
\ref{fig:ougammaoptkopt}a (green). The key point to look for is the rate at
which the information content is increased when clustering the analyzed data in
comparison to the rate when clustering the reference time series.

Figure \ref{fig:ougammaoptkopt}b shows the distance between the two clustering
curves.  A reasonable choice is $K_{\textrm{opt}}= 6$ because afterwards the
curve of the analyzed data does not approach the reference one,  which means no
significant further increase of the model diversity with respect to the
reflected Wiener process. The expected value is built based on $10$ reference
data sets. This number will be increased in the following examples to obtain
better estimates.

\begin{figure}[h]
\centering
\includegraphics[width=1.0\linewidth]{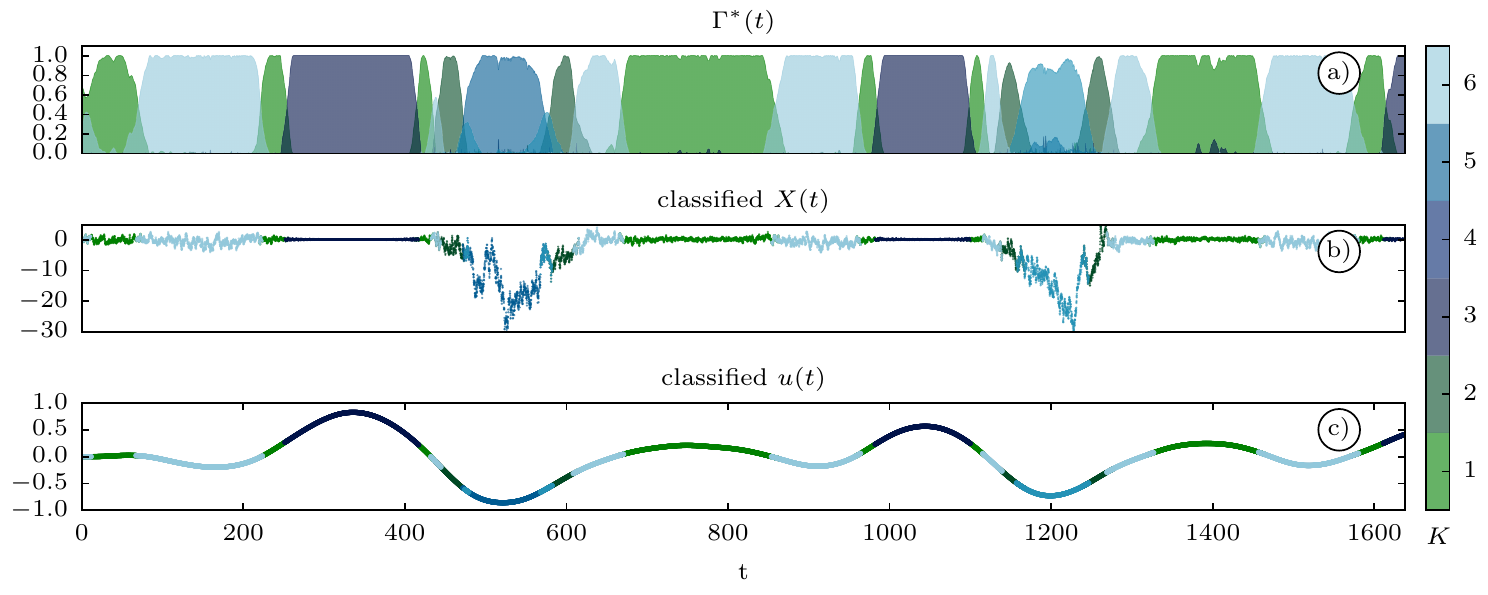}
\caption{The result of clustering the example time series that is presented in 
Fig. 
\ref{fig:OU_example}. a) The affiliation function that is found in solving the 
variational problem Eq. \eqref{eq.:FEM_discetisation} with $K=6$ and 
$\epsilon^2 \approx 20$. b) The classified time series $X(t)$. c) The classified 
auxiliary process $u(t)$. }
\label{fig:OU_example_result}
\end{figure}

Figure \ref{fig:OU_example_result}a displays $\bm{\Gamma}^{\ast}(t)$ obtained
after clustering the time series $X(t)$ with $K=6$ (see Fig.
\ref{fig:OU_example}a). The result is considered to be valid, since at almost every time instance we observe a high value of $\bm{\Gamma}^{\ast}(t)$ indicating a high certainty of the parameters. In Figure
\ref{fig:OU_example_result}b the affiliation function is used to label the
training data. Note, although the sample path of $X(t)$ at $t\approx500$ appears
to look similar to the sample path at $t\approx1200$ the methodology is capable
to differentiate a change in the parameters of the SDE and separate these times
in two different clusters. This is attributed to the high accuracy of the
method.

In Figure \ref{fig:OU_example_result}c the affiliation function is used to
classify the auxiliary process $u(t)$. The regime occupation duration is
changing across clusters. One can note that it is related to the rate of change
of the auxiliary process. At times when this process is changing fast the
occupation time is short and when the process is almost steady the regime
occupation time is long. This observation highlights that the method is able to
tract relatively rapid changes in the parameters if they occur frequently.

\begin{figure}[h]
\centering
\includegraphics[width=1.0\linewidth]{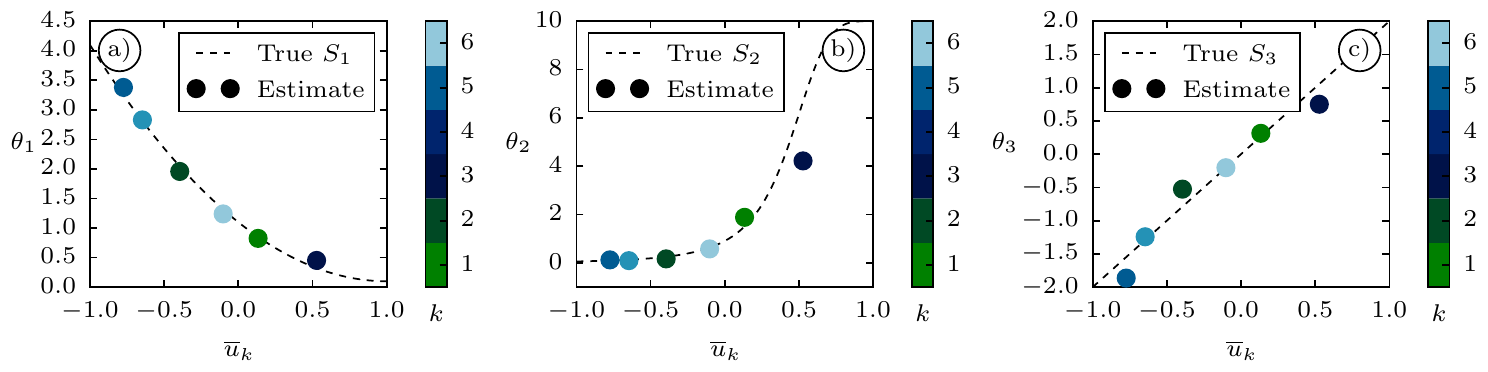}
\caption{The estimated parameters plotted over the cluster averaged auxiliary 
process. The 
model identification is performed with an optimal number of clusters equal to 
six. The 
colorbar is labeling the different clusters as in Figure 
\ref{fig:OU_example_result}. The dots represent the estimated 
parameter values. The dashed line shows the true scaling functions, which were 
used in generating the training data.}
\label{fig:OU_parameter_result}
\end{figure}
Figure \ref{fig:OU_parameter_result} shows the estimated parameters values for
six clusters. The corresponding six points are sufficient to estimate the hidden
scaling functions. When working with unexplored data the scaling functions would
need to be parameterized from the scatter plots by applying a further regression
analysis. We skip this step here because the estimated parameters are in a good
agreement with the true functions. Parameter $\theta_2$ shows largest
discrepancy in cluster $k=3$ (the right most point in fig.
\ref{fig:OU_parameter_result}b) which is explained by the fact, that this
cluster shows the largest variability of the process $u(t)$. Please observe the
classified periods of $u(t)$ for $k=3$ in figure \ref{fig:OU_example_result}c.
One way to improve the results would be to consider more sample points and
repeat the clustering with more clusters. \subsection{SDE with Two Independent Auxiliary Processes}
In this example, we consider an SDE with a nonlinear drift and a multiplicative
noise term. Furthermore, the goal of this example is to demonstrate that the
clustering methodology can handle modulation of two parameters which are
independent and have different timescales. The following functional form of 
the SDE is investigated:
\begin{equation}\label{eq:ND}
\textrm{d}X = (2 - \theta_1(t) X - \textrm{ln}(X^2)) \textrm{d}t +  
\theta_2(t) X \textrm{d}W\, ,  \qquad X(t_0) = 1 \, .
\end{equation}
The relationship between the parameters and the auxiliary processes $u(t)$ and
$v(t)$ is taken as
\begin{align}
	\theta_1(t) := S_1(u)&= (u(t)+1)^4+0.1\, ,\\
	\theta_2(t) := S_2(v)&= -2v(t) + 2\, ,
\end{align} 
where $u(t)$ and $v(t)$ are independent solutions of Eq. \eqref{eq:OU_U(t)} with
different values for $T_c$: $T_{c, v} = 4T_{c, u}$ (see Eq. \eqref{eq:OU_U(t)}).
The task is to recover the scaling functions solely from discrete-time
observations of $X(t), u(t)$ and $v(t)$.For details on running the optimization
algorithm for this example see Appendix C.
\begin{figure}[h]
	\centering
	\includegraphics[width=1.0\linewidth]{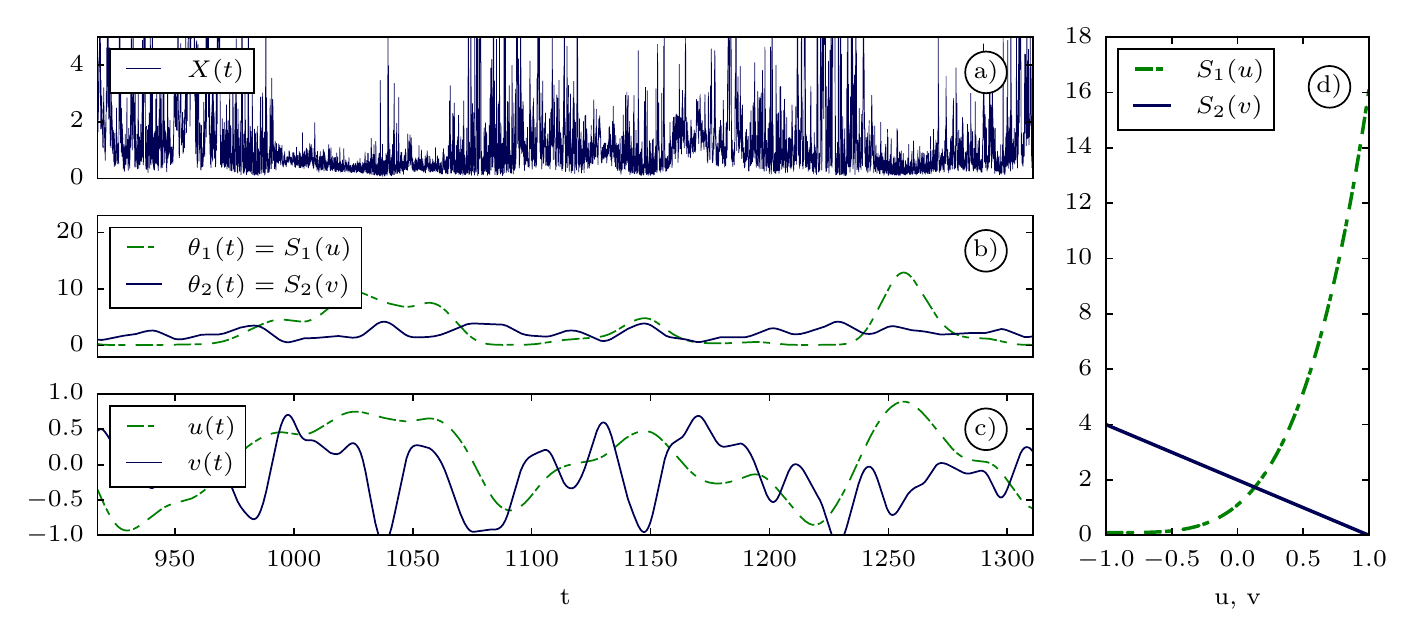}
	\caption{Summary of the functions to generate the non-stationary training 
	data 
	$X(t)$
		according to the Eq. \eqref{eq:ND}. a) The sample path of 
		$X(t)$. b) The
		temporal evolution of the model parameters. c) The auxiliary parameter 
		processes $u(t)$ and $v(t)$. d) The scaling functions $\theta_n(t) = 
		S_n(u(t))$.}
	\label{fig:ND_example}
\end{figure}

The sample size of the time series $X(t)$ is $131072$ ($\Delta t $ grid) points 
and it is longer
than in the example with the OU process (see Sec. \ref{sec.:ou_example}). More
points support the methodology to distinguish the parameter variability. The
processes $u(t)$ and $v(t)$ are uncorrelated and the re-occurrence of the unique
parameter pairs in the sample path $X(t)$ is not assured. This makes the
identification of the scaling functions ambitious. To mitigate the problem one
needs to increase the number of clusters to raise the resolution of the scaling
functions. By increasing the number of clusters we reduce the number of the
samples available to estimate the parameter values within one cluster.
Consequently, one also needs to increase the sample size of the total time series $X(t)$ to maintain the accuracy of the parameter estimates.
\begin{figure}[h]
	\centering
	\includegraphics[width=1.0\linewidth]{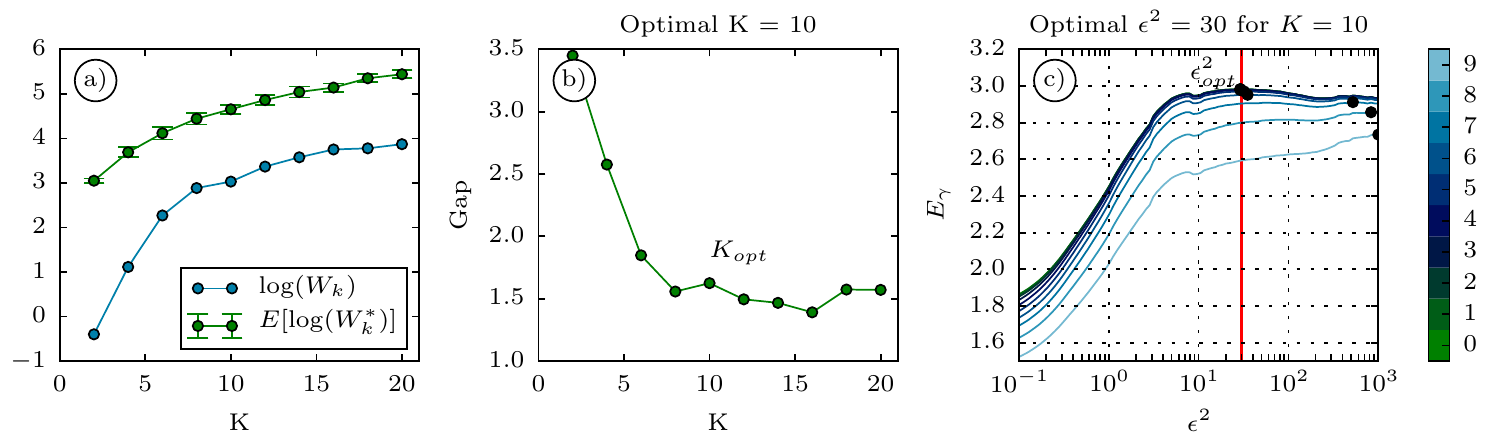}
	\caption{Estimation of the number of clusters for the example with two 
	auxiliary processes (see 
		Fig. \eqref{fig:ND_example} and Eq. \eqref{eq:ND}) using the gap statistics 
		approach a) (see Eq. \eqref{eq.:diversity_measure}) and b)(see Eq. 
		\eqref{eq.:gap}). In the
		panel a) after approximately $K=10$ the diversity of the 
		training data is 
		increasing at the same rate as in the reference data sets. This suggest 
		the 
		optimal $K$ to be around $10$.   
		c) Estimation of the regularization parameter $\epsilon^2$ for 
		clustering 
		with $K=10$. The colorbar codes the consecutive suppression of high 
		frequency bands to emphasize the  maximum. The filtering step is not 
		necessary 
		in this example.}
	\label{fig:ndgammaoptkopt}
\end{figure}

Figure \ref{fig:ND_example}a shows the sample path of the SDE \eqref{eq:ND}
together with the auxiliary processes $u(t)$ and $v(t)$ (see Fig.
\ref{fig:ND_example}c) which were used to generate $\bm{\Theta}(t) =
[\theta_1(t), \theta_2(t)]$ (see Fig. \ref{fig:ND_example}b). The scaling
functions are plotted in Fig. \ref{fig:ndgammaoptkopt}d.

The diversity measure of the clustering for the process $X(t)$ and the reference
data set is shown for comparison in Fig. \ref{fig:ndgammaoptkopt}a. The peak 
in diversity when clustering the data set $X(t)$ is reached
approximately at $K_{\textrm{opt}}=10$. By plotting the distance between the two
curves we find that for $K>10$ the clustering does not add significantly more 
 information relative to the clustering of the reflected Wiener process. The
selection of $\epsilon^2_{\textrm{opt}}$ for this example shows a flat maximum
(see Fig. \ref{fig:ndgammaoptkopt}c).
\begin{figure}[h]
	\centering
	\includegraphics[width=1.0\linewidth]{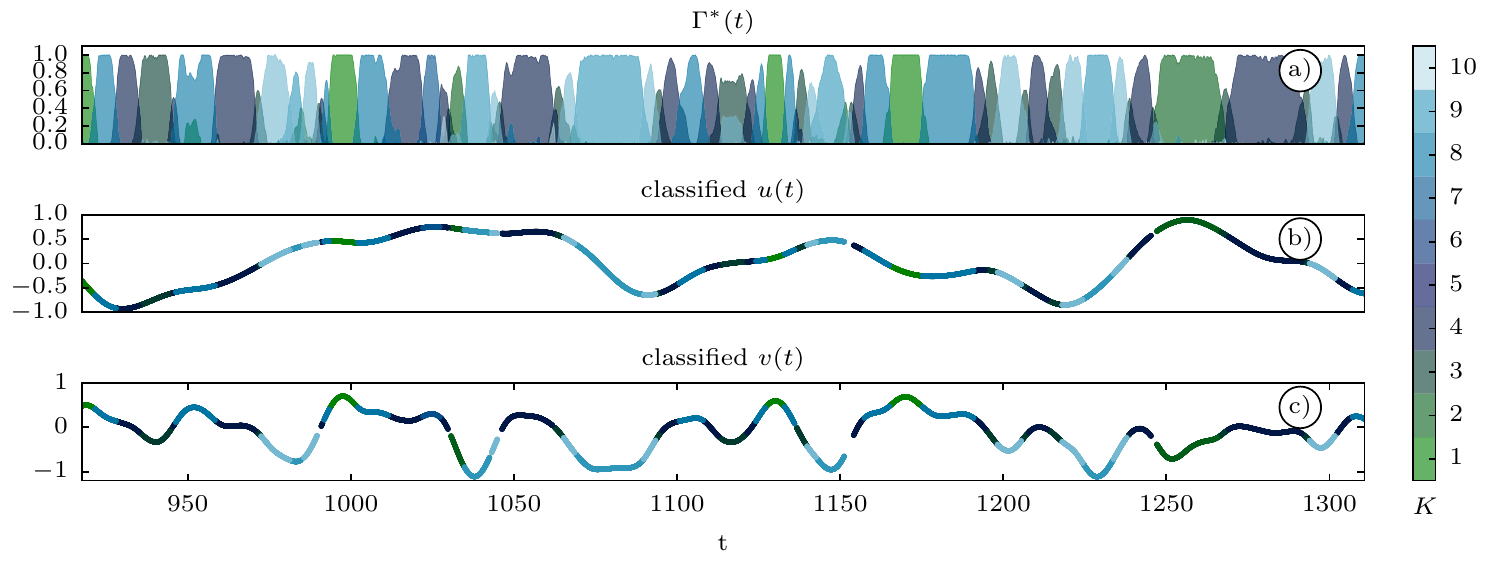}
	\caption{Result of the clustering the example time series that is 
	presented in 
	Fig. 
		\ref{fig:ND_example}. a) The affiliation function that is found in 
		solving the 
		variational problem \ref{eq.:FEM_discetisation} with $K=10$ and 
		$\epsilon^2=30$. b) The classified auxiliary process $u(t)$. c) The 
		classified 
		auxiliary process $v(t)$. The figure shows $30 \%$ of the 
		data set. }
	\label{fig:ndexampleresult}
\end{figure}

\begin{figure}[h]
	\centering
	\includegraphics[width=0.91\linewidth]{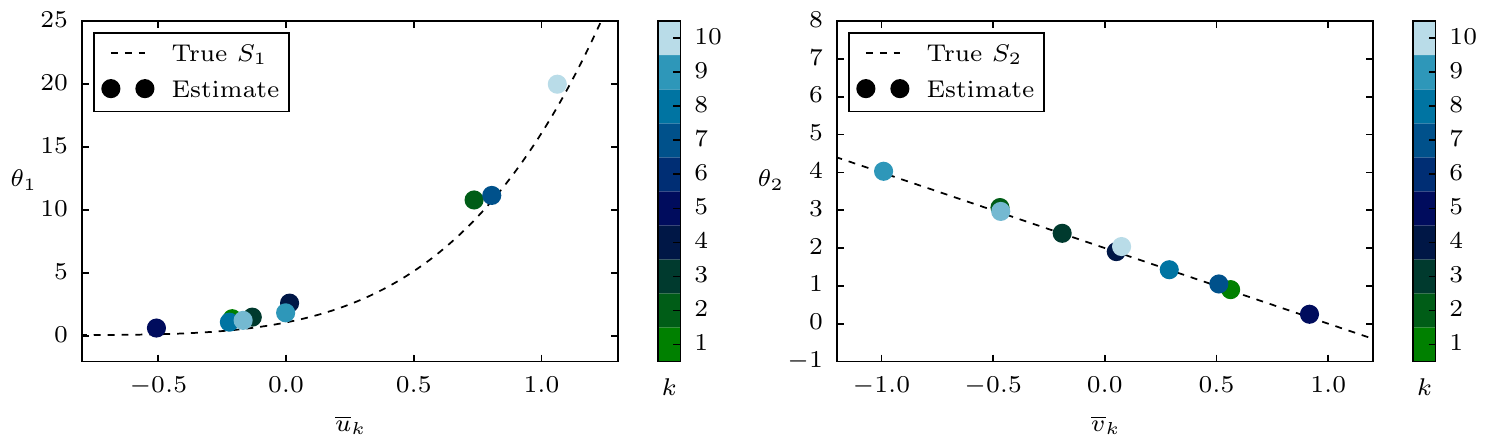}
	\caption{Estimated parameters plotted over the cluster averaged 
	auxiliary 
		processes. The 
		model identification is performed with an optimal number of clusters 
		equal to 
		ten. The 
		colorbar is labeling the different clusters as in Figure 
		\ref{fig:OU_example_result}. The dots represent the estimated 
		parameter values. The dashed line shows the true scaling functions, 
		which were 
		used in generating the training data.}
	\label{fig:ndparameterresult}
\end{figure}

\begin{figure}[h]
	\centering
	\includegraphics[width=1.0\linewidth]{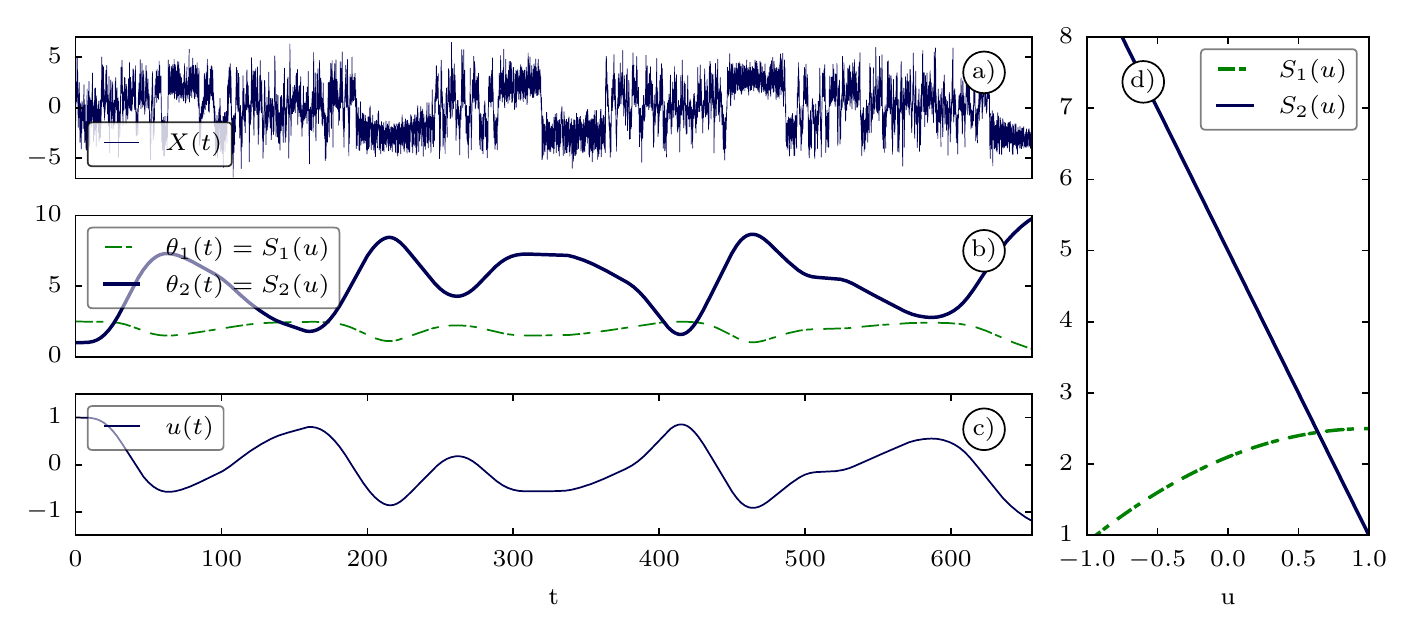}
	\caption{The example of a non-stationary process that is associated with Eq. 
	\eqref{eq:CN}. a) The sample 
		path 
		of $X(t)$. 
		b) The
		temporal evolution of the model parameters. c) The parameter auxiliary 
		process 
		$u(t)$. d) The scaling functions $ 
		S_n(u(t)) = \theta_n(t) $.}
	\label{fig:CN_example}
\end{figure}

The optimal classifier $\bm{\Gamma}^{\ast}(t)$ and the partitioning of the
auxiliary processes are shown in Fig.  \ref{fig:ndexampleresult} where we can
investigate the quality of the clustering. Figure \ref{fig:ndparameterresult}
displays the estimated parameters and the true scaling functions. By visually
inspecting these results we note a consistent spread of the estimates in the
parameter $\theta_2$, but estimates in the parameter $\theta_1$ are distributed
rather disproportionally. It is noticeable in the estimates of $\theta_2$  that
three parameters have been apparently estimated twice and are positioned close
to each other. The respective cluster values for the parameter $\theta_1$
however are located apart, such that the scaling function is sampled at
significantly different points. This means that in this example the methodology
is distinguishing pairs of parameters which have one common element. We
therefore conclude that the methodology is capable to classify SDE models with
at least two uncorrelated auxiliary processes.
In this example, we investigate an SDE with nonlinear drift and nonlinear
diffusion. For the sake of simplicity, we consider only one auxiliary
process $u(t)$ and two parameters: one in the drift term and one in the
diffusion term. The functional form of the SDE is
\begin{equation}\label{eq:CN}
\textrm{d}X = (\theta_1(t)X - X^3) \textrm{d}t + \theta_2(t)\sqrt{1 + X^2} 
\textrm{d}W\, , \qquad X(t_0) = 0\, .
\end{equation}
It is possible to incorporate more parameters into the model structure, for
example, in front the quadratic or the cubic term. The parameter $\theta_1(t)$
is a time-dependent bifurcation parameter which causes a continuous variation of
the equilibrium properties in time, eventually leading to metastable states in
the dynamics. (see fig. \ref{fig:CN_example}a). The clustering methodology can
comprehend transitions in the underlying double-well potential of the SDE. The
relationship between the parameters and the auxiliary processes is
\begin{align}
	\theta_1(t) \equiv S_1(u)&= -0.4(u(t)-1)^2 + 2.5\, , \\
	\theta_2(t) \equiv S_2(u)&= -4 u(t) + 5\, .
\end{align} 
Like in the previous examples the task is to recover the scaling functions by
only knowing one discrete-time trajectory of the process $X(t)$ and $u(t)$.

\begin{figure}[h]
	\centering
	\includegraphics[width=1.0\linewidth]{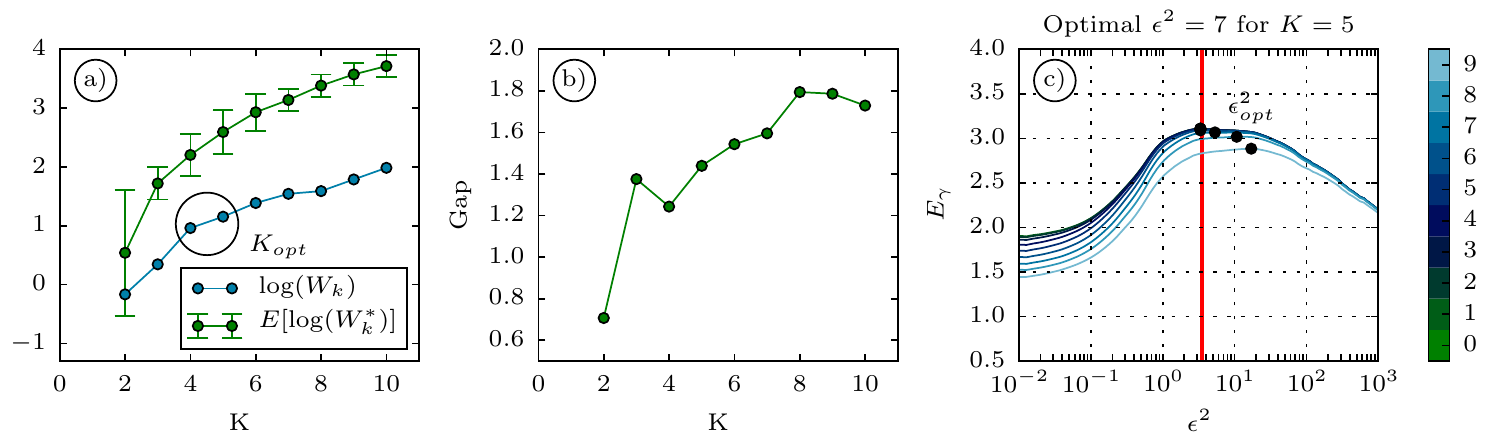}
	\caption{An estimation of the number of clusters 
	for the example associated with Eq. \eqref{eq:CN} (see 
		Fig. \ref{fig:CN_example}) using the gap statistics 
		approach a) (see Eq. \eqref{eq.:diversity_measure}) and b)(see Eq. 
		\eqref{eq.:gap}). In 
		panel a) we find that for $K > 4$ the diversity of the training data is 
		increasing at constant rate. b) The gap statistics suggests 
		$K_{opt}=4$.
		c) The estimation of the regularization parameter $\epsilon^2$ for 
		clustering 
		with $K=5$. The colorbar codes the consecutive suppression of high 
		frequency bands to emphasize the  maximum.}
	\label{fig:cngammaoptkopt}
\end{figure}

In this example our method to estimate $K_{\textrm{opt}}$ does not provide such
a clear picture as it did in the previous examples. Nevertheless, it is still
beneficial to analyze the tendencies of the curves in Figure
\ref{fig:cngammaoptkopt}a, b. As indicated by the minimum in the Figure
\ref{fig:cngammaoptkopt}b we should choose $K_{\textrm{opt}}=2$. However, our
final goal is to recover the scaling functions and in the case of a nonlinear
scaling functions one requires more than two points and therefore more
clusters. According to the second minimum in Figure \ref{fig:cngammaoptkopt}b,
the second-best choice is $K_{\textrm{opt}}=4$. If we investigate the curve
$\textrm{log}(W_k)$ in Figure \ref{fig:cngammaoptkopt}a, we find a minor
characteristic change in the slope at $K_{\textrm{opt}}=4$. To resolve the
scaling functions even better we select $K_{\textrm{opt}}=5$.

\begin{figure}[h]
	\centering
	\includegraphics[width=1.0\linewidth]{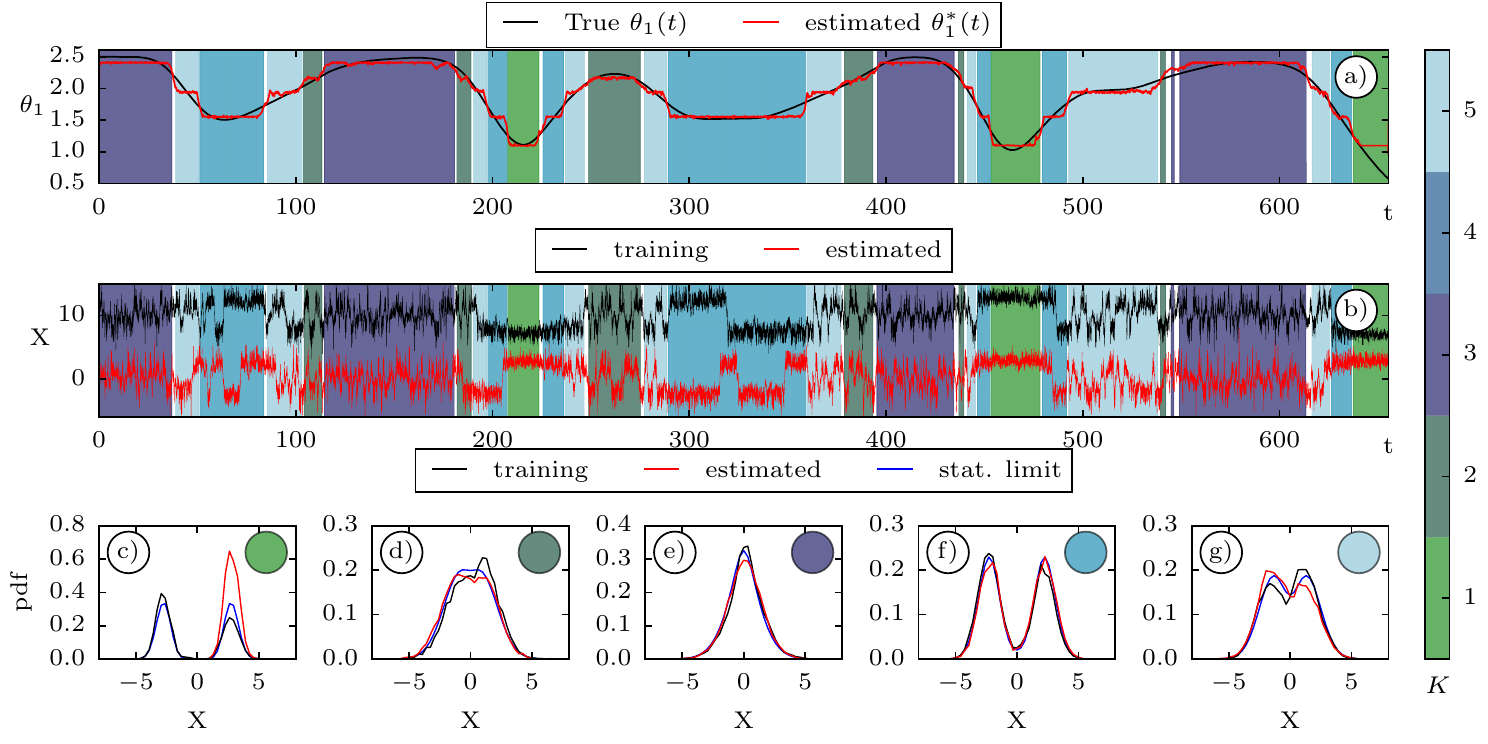}
	\caption{Clustered time series generated with 
	Eq. \eqref{eq:CN}. In a) and b) the colored background is the estimated 
	affiliation vector $\bm{\Gamma}^{\ast}(t) > 0.8$. The colorbar to the right of the 
	panels labels the corresponding $5$ clusters. a) The true temporal 
	evolution (black) of the drift parameter and its recovered version 
	(red). Other values of the hyperparameters are given in Fig. \ref{fig:cnexampleresultmultytheta}. The second parameter looks similar and is not shown. b)The 
	simulation 
	of 
	the sample path with the estimated parameter values compared to the 
	training data. Panels c) to g) show the cluster respective densities: the 
	histograms for the training 
	(black) and simulation (red). The blue lines mark the stationary pdf which 
	is computed from 
	the parameter estimates. The colored circles are 
	marking the corresponding clusters.}
	\label{fig:cnexampleresult}
\end{figure}
Knowing the affiliation function and the averaged
values of the parameters,  the temporal evolution of the model parameters (see
Fig. \ref{fig:cnexampleresult}a) is approximated by
\begin{equation}\label{eq.theta approx}
	\theta^{\ast}(t) = \sum_{k=1}^{K_{\textrm{opt}}}\overline{\theta}^{\ast}_k 
	\gamma_k(t)\, ,
\end{equation}
where $\overline{\theta}^{\ast}_k$ is a constant vector of size $N$. Figure
\ref{fig:cnexampleresult}a shows that the methodology recovers the hidden
evolution of the parameters rather well. The approximation becomes better the
more clusters are used (not shown). The drawback is that the uncertainty in
$\overline{\bm{\Theta}}^{\ast}$ grows with the number of clusters, due to the 
reduction in data points per cluster.

The approximated path of the model parameters (see Fig.
\ref{fig:cnexampleresult}a in red) is used to test the prediction performance
(see Fig. \ref{fig:cnexampleresult}b in red). In this figure, the time series of
the training dataset is shifted up for better comparison by the constant value
of $10$. The considered SDE exhibits a temporal modulation by the bifurcation
parameter. The times with two metastable states are expected to have different
sample paths because the local state of the system is dependent on the
particular realization of the Wiener process. Two different Wiener processes
were used to generate the training and simulation solution. Figure
\ref{fig:cnexampleresult}b shows the properly captured signal dynamics. The
corresponding cluster probability density functions are in good agreement with
those of the training data and with the respective stationary distributions,
which were computed from the parameter estimates. The largest discrepancy is
shown for the cluster with index $1$ (see Fig. \ref{fig:cnexampleresult}c). This
cluster has a short lifetime to capture the regime jumping in the present
realization. The trajectory-based p.d.f. estimate is off due to not observing
the entire cluster dynamics. If we consider the stationary distribution computed
from the estimated parameter values, we find the correct bimodal distribution
(compare in Fig. \ref{fig:cnexampleresult}c).

\begin{figure}[h]
	\centering
	\includegraphics[width=1.0\linewidth]{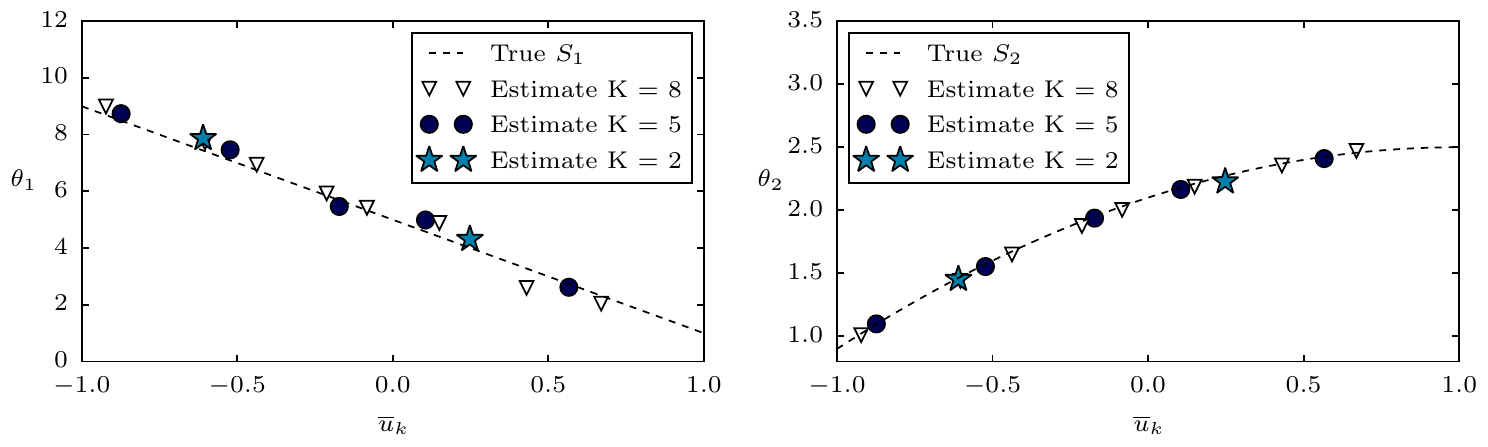}
	\caption{The estimated parameters plotted over the cluster averaged 
	auxiliary 
		process. The 
		model identification is performed with an optimal number of clusters 
		equal to 
		six. The dots represent the estimated 
		parameter values. The dashed line shows the true scaling functions, 
		which were 
		used in generating the training data. The optimal number of clusters is of second importance if one is able to parameterize the scaling functions.}
	\label{fig:CN_parameter_result}
\end{figure}

Figure \ref{fig:cnexampleresult} validates the performance of the model within
the training data. To construct a self-contained, non-stationary prediction
model, we show that the scaling functions between the auxiliary process $u(t)$
and the parameters are recovered as well (see Fig.
\ref{fig:CN_parameter_result}).
The accuracy of the clustering method is sufficient to obtain the correct
parameter values for a wide range of clusters used. It is important to note that
the optimal number of clusters is of secondary importance. The number of
clusters should be as low as possible and high enough to identify the scaling
function. After the relation between the auxiliary process and the model 
parameter is identified ( see Fig. \ref{fig:CN_parameter_result}), the
last step is consisting of performing a regression analysis to parameterize the
scaling behavior. It is left out here as it will be specific for each
application.

 \clearpage
  
\section{Discussion and Conclusion}
Data-driven modeling has the ability to enhance the understanding of complex
dynamical systems by highlighting patterns or hidden regularity. Especially
challenging is the identification and the parameterization of nonlinear and
non-stationary SDE. The intention of the system identification step is the
determination of a suitable model structure and the estimation of the
corresponding  model parameters. In retrospect, one may analyze the data
associated with the newly estimated model to improve understanding of the
process to derive hidden relationships. These relationships may not be given
explicitly by the clustering approach, so that unveiling these hidden
dependencies requires additional research work. For instance,
\cite{vercauteren_investigation_2016} characterized the influence of
non-turbulent motions on the turbulence after a model-based data-clustering
methodology was applied \citep{vercauteren_clustering_2015}. In a study of the
North Atlantic Oscillation in the climate system, \cite{quinn_dynamical_2021}
used the clustering results to study dynamical processes associated with
atmospheric regime transitions.

The presented approach clusters the data into segments and models each of them
with an individual, locally-stationary model. Such a decomposition of a complex
dynamical behavior may not be optimal since it may not always be clear how to
differentiate between a nonlinear and non-stationary time series without
knowing the underlying structure of the data-generating process. Identifying an
appropriate model structure is therefore an important step. For discrete-time
systems, \cite{billings_nonlinear_2013} developed an estimation algorithm which
can sequentially rate and collect important model terms. Following that,
\cite{wei_term_2004} find an effective model for a highly nonlinear
terrestrial magnetospheric dynamical system. Within their framework the
time-varying parameters can be estimated with a multi-wavelet
\citep{wei_time-varying_2010} or a sliding window approach
\citep{li_identification_2016}.

In section 2 we recall the variational clustering framework which allows the
simultaneous analysis and modeling of a non-stationary time series. One
variation of the framework \citep{horenko_identification_2010} is based on
discretized mathematical models, which are also termed discrete-time systems. To
integrate an identified discrete-time model into a preexisting continuous-time
system, one needs to account for the discretization operation, otherwise, the
solution is restricted to the time frequency given by the data. One way to
overcome the problem is to convey the frequency response of the identified model
by transforming it into a continuous one \citep{kuznetsov_simple_1994}. This
transformation is widely studied in linear system theory and is extendable for
the case of nonlinear systems \citep{billings_nonlinear_2013}[p. 342]. Our
approach avoids the above-mentioned difficulty because we use the MLE derived
from an appropriate likelihood (fitness) function for the discrete-time
observations of SDE.

In general, the  main challenge in applying the MLE is the missing transition
density function, which requires, for example, a solution of the Fokker-Plank 
equation for a
considered structure of the SDE. Especially for the nonlinear SDE, this function
is analytically not known and demands estimation \citep{fuchs_inference_2013}. 
In our approach, by separating the available time series into $K$ clusters,
the available number of samples for estimation reduces, leading thereby to
uncertain estimates in each of them. The closed-form approximation based on the
Hermite expansion provides an accurate approximation for the transition density,
with an error of $O(\Delta t^4)$ \citep{ait-sahalia_maximum_2002}. This
accuracy allows to increase the number of clusters and to raise the resolution
of the hidden scaling functions making an accurate parameterization task
feasible.

We combine the closed-form likelihood function approach based on suitable
Hermite-expansions with the non-parametric clustering framework. Specifically,
we highlight the details of the necessary modification in the original subspace
clustering algorithm to achieve accurate estimates of the model parameters
independently of the discretization of the QP problem. Furthermore, we present an
extensive numerical study consisting of three controlled examples of different
complexity to validate both the parameterization methodology and the novel
methods to estimate the hyperparameters of the framework, which are independent
and does not rely on information theory.

In particular, our approach to estimating the optimal number of clusters is
based on measuring the degree of model diversity and comparing it to the
clustering of the reflected reference Wiener process. The numerical experiments
confirm that this approach produces adequately reliable estimates. However, the
associated graphs prove to be dependent on the considered model structure such
that the selection requires a case-specific interpretation. The proposed method
is still computationally expensive because it requires to repeat the clustering 
of the reference
data set multiple times to construct an expected value.

Our method to identify the optimal regularization value of the affiliation
function proves to be robust for the studied test cases. We further extended it
to articulate the optimum through consecutively filtering out the small-scale
variability from the affiliation vector. Unfortunately, we could not construct an appropriated example to show its potential. However, apart from this study, we
tested this methodology on the real-data application and found it to be effective. This will be reported in a follow-up study.

In summary, the data-driven variational clustering approach enables modeling of
the nonlinear and non-stationary processes to develop models for multi-scale
dynamical systems. The model structure of the SDE is thereby unknown and demands
to be prescribed or inferred in supplementary ways. The parameterization method
is accurate and provides reasonable results with small data sets consisting of
$2\times10^4$ points. It is also scalable on the GPU as it is realized in this
work. We tested data sets consisting of $2\times10^6$ points, where the
computation time was approximately 24 hours on a GeForce 1080 with 3 clusters.
The numerical examples demonstrate the correct identification of the required
scaling functions, which allows parameterizing a considered SDE by linking
additional auxiliary processes to the model parameters.

This work focuses on one-dimensional time series; however, it is possible to extend the approach to a multi-dimensional case. One should consider following two research directions. The first one is to use the closed-form likelihood expansion for multivariate diffusion \citep{ait-sahalia_closed-form_2008}. Modification of the variational clustering method is not required for this. However, practically one may run into problems of over-parameterization because the cross-dimensional interaction terms and their temporal variability lead to unconstrained systems. One way to mitigate this is to enforce diagonally-spars diffusion and drift matrices as it is demonstrated by \cite{okane_research_2017-1}. 

The second research direction is to examine multi-dimensional or even spatial problems. Then it is suggestible to replace the MLE with an ensemble Kalman inversion method \citep{albers_ensemble_2019}. The variational clustering framework is also generalized to classify the spatial parameter variability \cite[p. 81]{kaiser_data-based_2015}. \cite{pospisil_scalable_2018} implemented this approach within the openly available code, which is linked in their work. The proposed methods appear promising, however with increasing dimensions, it becomes more and more important to estimate the correct functional terms rather than the parameter values. Hence, an additional inverse problem may be combined with the presented approach to assessing the structural form of the system.

Finally, it is noteworthy that the presented parameterization method is intended
to be used to develop a data-driven stochastic parameterization of intermittent
turbulence present in the nocturnal atmospheric boundary layer. In climate
models for example, the limited numerical resolution leaves unresolved degrees
of freedom, commonly denoted as subgrid-scale processes. The approach presented
here provides an ideal framework to derive subgrid-scale models which are
modulated by resolved variables. These results will be published elsewhere.
Another direction for future research is the extension of the non-stationary data-driven learning approach developed in this work to SDE models with (at least) two widely separated time scales, for which the MLE is known to become asymptotically biased. For these systems the proposed methodology has to be combined with multiscale robust inference techniques
\citep{krumscheid_semiparametric_2013, kalliadasis_new_2015, krumscheid_perturbation-based_2018}.

 \clearpage
\appendix
\section{Reusing of the initial guess}
One can estimate the values of the hyperparameters $K$ and $\epsilon^2$ by
running the optimization problem for each pair of parameter combinations. In the
OU example, we choose for $K $ the range of values $[ 2,3, \dots, 10]$ and for
$\epsilon^2 $ the range of values $[10^{-1}, 10^2 ]$, where we use $100$
divisions on the logarithmic scale for the latter one. In principle for each
parameter combination $(K, \epsilon^2)$ the optimization method can be started
with a random initial guess $\overline{\bm{\Theta}}_0$ and $\bm{\Gamma}_0(t)$.
However, practically it is efficient to use the optimal values from the previous
run in the following way:
\begin{enumerate}
	\item Select one value for $K$,
	\item Initialize the first optimization run with the smallest or the 
	largest 
	value 
	of 
	$\epsilon^2_0$ together with a random guess for 
	$\overline{\bm{\Theta}}_0$ and 
	$\bm{\Gamma}_0(t)$,
	\item Algorithm 1 provides a solution 
	$\overline{\bm{\Theta}}^{\ast}_1$ 
	and 
	$\bm{\Gamma}^{\ast}_1$
	\item Change the parameter value $\epsilon^2$ to the next closest value 
	$\epsilon^2_1$,
	\item Start the optimization method with $\epsilon^2_1$ and the initial 
	guess to be 
	the solution from the previous run: $\overline{\bm{\Theta}}^{\ast}_1$ 
	and $\bm{\Gamma}^{\ast}_1$,
	\item Algorithm 1 provides to the solution 
	$\overline{\bm{\Theta}}^{\ast}_2$ and $\bm{\Gamma}^{\ast}_2$,
	\item Solve the variational problem for all considered $\epsilon^2$ by 
	incrementally changing its value,
	\item Change the value of $K$ and repeat with (1).
\end{enumerate}
The approach where the initial guess is taken randomly requires , of course,
more computational time in total than the case with reused solution of previous
simulations. The principle of traversing the regularization parameter
$\epsilon^2$ is further exploited by traversing the parameter from low to high
value and back to a low value while consistently reusing the previous solution 
outcomes.
This strategy showed to be effective to find a better solution after each
iteration, although with only minor incremental improvement.

\section{Details on the settings of the framework for the OU example in section 
5.3}
The minimization algorithm for the QP solver has two stopping criteria which are
as follows. The maximum number of iterations is set to $100$, and the minimum
difference between consecutive cost function values is set to $10^{-8}$. The
reduction value for the QP solver speeds up the QP problem and is set to $1/3$.
For the $\bm{\Theta}$ solver, the global optimizer is a random search algorithm
supplemented with a local search algorithm. For both of them, the maximum
number of function evaluations is set to $300$ and the break tolerance to a 
relative value 
$10^{-10}$. The
population size of the global optimizer is set to $3000$. The final point to
mention is the  parameter bounds for the $\bm{\Theta}$ solver. In this example (and for
the clustering of the reference data set) they are :$-20<\theta_1<20;\
0<\theta_2<20;\ 0<\theta_3<20$.

\section{Details on the settings of the framework for the example in Section 
5.4}
Here are the details on the simulation of the sample path $X(t)$ (see Fig.
\ref{fig:CN_example}) which forms the training data for the clustering test. The
time series $X(t)$ has $65536$ samples and a constant time step $\Delta t =
0.01$. The sample path  is obtained by solving the Eq. \ref{eq:CN} with the
Milstein method and the sample path of $u(t)$ is obtained by solving the Eq.
\ref{eq:OU_U(t)} with the Euler-Maruyama method. During the simulation the time
step between each of the $65536$ samples is reduced to the the value $10^{-4}$
in order to obtain a more accurate realization of the sample path. The initial
value for Eq. \ref{eq:CN} is set to $X(t_0) = 0$ and for Eq. \ref{eq:OU_U(t)} -
to $U(t_0) = [1, 0, 0, 0]$. The parameters in Eq. \ref{eq:OU_U(t)} were taken 
to be $T_c = 15$; $b_0 = 0.2$.

Next we provide the details on parameters of the clustering framework itself. 
In this example, the maximum number of iteration is set to $100$ and the minimum
difference in the consecutive evaluation of the fitness function is set to a
relative value $10^{-8}$. The reduction value for the $\bm{\Gamma}$ solver is set to
$\alpha = 0.1$. For both of solvers, the maximum number of calls is set to $500$
and the termination tolerance to a relative value $10^{-10}$. The population
size of the global optimizer is set to $1000$. The bounds for the parameter
limits of the $\bm{\Theta}$ solver (and for the clustering of the reference data set)
is:$-10<\theta_1<10;\ 0<\theta_2<10$.

\section{Coefficients of the  Hermit Expansion }
The coefficients $\eta_j$ in Eq. \eqref{eq.:expantion} are approximated with a
truncated Taylor series of length $M = J/2$ , where $J = 6$. We state them here 
for completeness.
\begin{equation*}
\begin{split}
\eta_{0}^{(3)} &= 1\;,\\
\eta_{1}^{(3)} &= -\mu {h}^{1/2} - \frac{ 2\mu\mu' + \mu''}{4}
{h}^{3/2} - \frac{ 4\mu {\mu'}^{2} + 4 {\mu}^{2}\mu''+ 6\mu'
	\mu''+4
	\mu \mu^{(3)}+ \mu^{(4)}}{24} {h}^{5/2}\;,\\
\eta_{2}^{(3)} &= \frac{{\mu}^{2}+\mu'}{2} h + \frac{ 6{\mu}^{2}\mu'+4
	{(\mu')}^{2}+7\mu \mu''+2\mu^{(3)}}{12} {h}^{2} +
\frac{1}{96} \Bigl( 28{\mu}^{2}{(\mu')}^{2}+28{\mu}^{2} \mu^{(3)}\\
&\quad + 16{(\mu')}^{3}+16{\mu}^{3}{\mu''}+88 \mu \mu'
\mu''+21{(\mu'')}^{2}+32\mu' \mu^{(3)}+16 \mu
\mu^{(4)}+3\mu^{(5)}\Bigr){h}^{3}\;,\\
\eta_{3}^{(3)} &= - \frac{{\mu}^{3}+3\mu \mu'+\mu''}{6} {h}^{3/2} -
\frac{1}{48}\Bigl( 12{\mu}^{3}\mu'+
28\mu{(\mu')}^{2}+22{\mu}^{2}\mu''+24
\mu'\mu''\\
&\quad +14 \mu \mu^{(3)}+3 \mu^{(4)}\Bigr) {h}^{5/2}\;,\\
\eta_{4}^{(3)} &= \frac{ {\mu}^{4}+6{\mu}^{2}\mu'+3{(\mu')}^{2} +4 \mu \mu''
	+ \mu^{(3)} }{24} {h}^{2} + \frac{1}{240}\Bigl( 20
{\mu}^{4}\mu' +50
{\mu }^{3}\mu'' + 100 {\mu}^{2} {(\mu')}^{2}\\
&\quad +50 {\mu}^{2} \mu^{(3)} + 23 \mu \mu^{(4)} + 180 \mu \mu'
\mu'' +40 {(\mu')}^{3} + 34 {(\mu'')}^{2} + 52 \mu' \mu^{(3)} + 4
\mu^{(5)} \Bigr) {h}^{3}\;,\\
\eta_{5}^{(3)} &=
-\frac{{\mu}^{5}+10{\mu}^{3}\mu'+15\mu{(\mu')}^{2}+10{\mu}^{2}\mu''+
	10\mu' \mu''+5\mu \mu^{(3)}+\mu^{(4)} }{120} {h}^{5/2}\;,\\
\eta_{6}^{(3)} &= \frac{1}{720}
\Bigl({\mu}^{6}+15{\mu}^{4}\mu'+15{(\mu')}^{3}+20{\mu}^{3}\mu''+15\mu'\mu^{(3)}+45{\mu}^{2}{(\mu')}^{2}+10{(\mu'')}^{2}\\
&\quad +15{\mu}^{2}\mu^{(3)}+60\mu \mu' \mu'' + 6 \mu \mu^{(4)}
+\mu^{(5)}\Bigr) {h}^{3}\;,
\end{split}
\end{equation*}
where $\mu$ is the drift of the $Y$ process evaluated at $y_0$ (see Eq. 
\eqref{eq.:fix_condition}), ''$^\prime$'' denotes a derivative and $h = \Delta 
t$. 
The Hermite polynomials for our expansion are:
$H_0(z) = 1, H_1(z) = - z, H_2(z) = -1 + z^2, H_3(z) = 3z-z^3, H_4(z) = 3-6z^2+z^4, H_5(z) = -15z + 10z^3 - z^5, H_6(z) = -15 + 45z^2 -15z^4 + z^6$.

\section{Variation of the Hyperparameters }
The section contains a sensitivity analysis to demonstrate different solution options of the clustering frame work.  
\begin{figure}
	\centering
	\includegraphics[width=1.\linewidth]{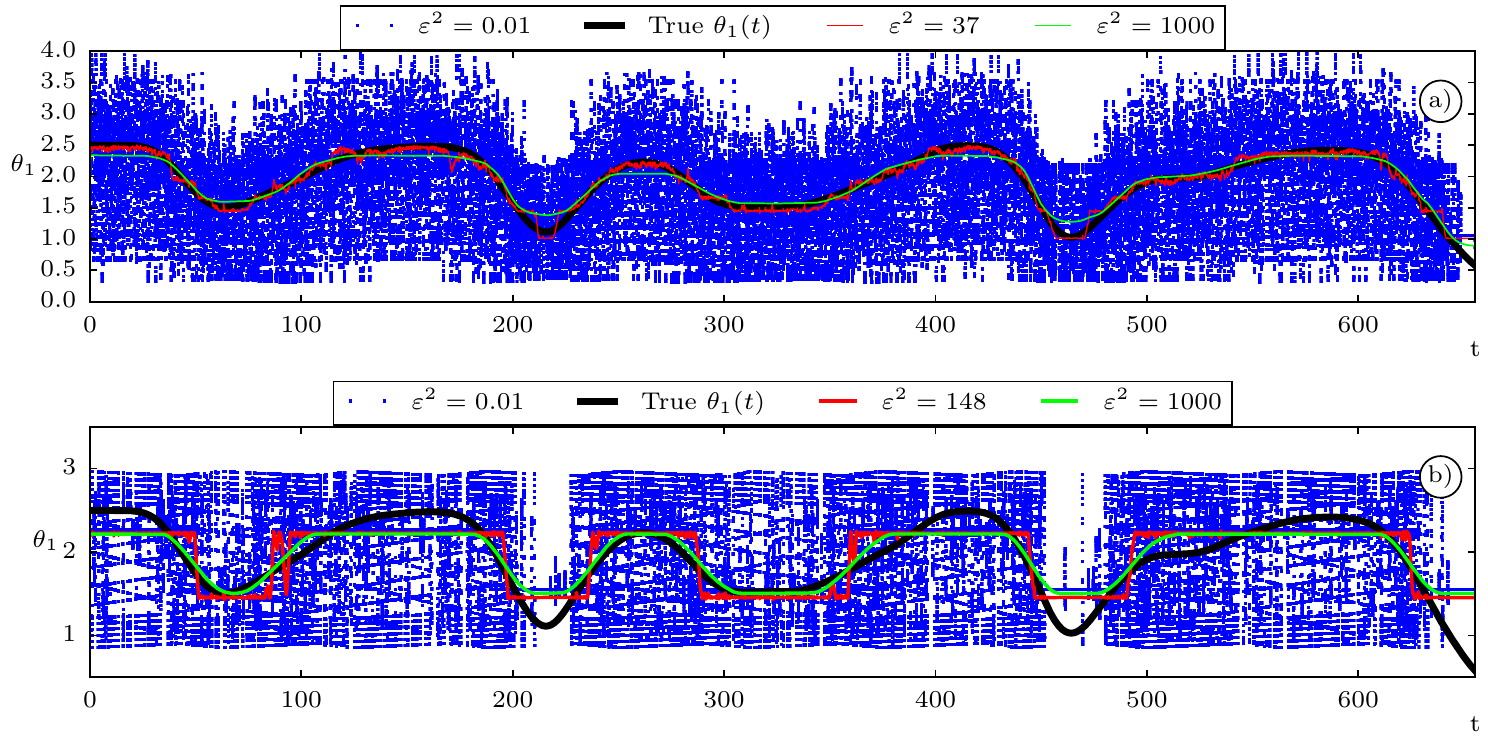}
	\caption{Reconstruction of the hidden evolution of one parameter for the nonlinear SDE using different values of the two significant hyperparameters ($K$ and $\varepsilon^2$). Panel a) shows the parameter reconstructed with eight clusters and panel b) using two clusters. Highly scattered blue dots label the weakly regularized solution. The red curves are suggested as an optimal value following the presented approach.}
	\label{fig:cnexampleresultmultytheta}
\end{figure}
 
\bibliographystyle{apalike} 

\bibliography{library2}

\end{document}